%% file: main.tex
\documentclass[review,onefignum,onetabnum]{siamart171218}
\usepackage{float}
\usepackage{amsmath,amsfonts,amssymb,mathtools,bm}
\usepackage{booktabs}
\usepackage{graphicx}
\usepackage{enumitem}
\usepackage{url}
\usepackage{xcolor}
\usepackage{subcaption}
\usepackage{array}

\hypersetup{
    colorlinks=true,
    linkcolor=blue,
    citecolor=blue,
    urlcolor=blue
}

\graphicspath{{figures/}{manuscript/figures/}{results/figures/}}

\newsiamremark{assumption}{Assumption}
\newsiamremark{example}{Example}

\colorlet{ORANGE}{orange}
\colorlet{GREEN}{green}
\colorlet{BLACK}{black}

\headers{Observable-Reduction-Guided Sparse Regression}{K. C. Nguyen and K. B. Flores}

\newcommand{\udot}{\dot u}
\newcommand{\uddot}{\ddot u}
\newcommand{\vdot}{\dot v}
\newcommand{\vddot}{\ddot v}
\newcommand{\gamone}{\gamma_1}
\newcommand{\gamtwo}{\gamma_2}

\newcommand{\ThetaLib}{\Theta}

\title{Observable-Reduction-Guided Sparse Regression for Partially Observed Active-Quiescent Systems}

\author{
Kyle C. Nguyen\thanks{Center for Research in Scientific Computation and Department of Mathematics, North Carolina State University, Raleigh, North Carolina, USA. Corresponding author: kcnguye2@ncsu.edu.}
\and
Kevin B. Flores\thanks{Center for Research in Scientific Computation and Department of Mathematics, North Carolina State University, Raleigh, North Carolina, USA. kbflores@ncsu.edu.}
}

\begin{document}
\nolinenumbers

\maketitle

\begin{abstract}
Active-quiescent switching occurs in biological populations where growth acts through a proliferative active state while cells may reversibly enter a nonproliferative quiescent state. Experiments often observe only part of this process, including direct or sparsely sampled active-state marker observations, aggregate-only population measurements, and aggregate measurements supplemented by sparse active-state marker data. These cases arise, for example, when only a limited marker panel is available, when active-state fluorescence requires fixation or endpoint assays such as flow cytometry, or when experiments report only total cell number, optical density, tumor burden, or aggregate fluorescence. Such observation choices complicate sparse regression-based equation learning methods, including sparse identification of nonlinear dynamics (SINDy), because the measured variable need not satisfy the same equation as the hidden active-quiescent system. A library chosen for the wrong observable may therefore fit a trajectory while losing mechanistic interpretability or transferability.

We study this issue in a two-compartment ordinary differential equation model. We call the elimination of hidden states observable reduction, meaning the differential equation satisfied by the measured variable after unobserved compartments have been removed. We derive these reductions for several biologically relevant and commonly used functional forms for growth, and use them to construct observation-specific sparse-regression libraries. Using computationally generated data, we compare these structured libraries with SINDy in its common polynomial-library form. In the examples studied here, polynomial libraries can match training trajectories while failing coefficient-relation and transfer checks. Libraries derived from the active or aggregate observable reductions recover interpretable coefficient maps when the observed regime is informative. These results show that interpretable equation learning for hidden-compartment systems requires matching the regression target and candidate library to the observation process.

\medskip
\noindent\textbf{Relevance to Life Sciences.}
Active-quiescent switching is common in biological populations whose cells, organisms, or subpopulations move between proliferative and nonproliferative states, including dormant, resting, or phenotypically inactive states. Experiments often do not observe all compartments directly. Instead, they may measure active-state markers, total population counts, optical density, tumor burden, or aggregate fluorescence. This work shows how such partial or aggregate observations change the equation that should be learned from data and develops an observable-reduction workflow for interpretable sparse regression in quiescent biological models.

\medskip
\noindent\textbf{Mathematical Content.}
The paper studies sparse-regression library construction for partially and aggregately observed active-quiescent ODE systems. For each observation map, we derive the scalar observable reduction that determines the appropriate regression target and candidate library. Active-state observations lead to second-order observable equations with coefficient-to-parameter relations, while aggregate observations lead to inverse-dependent or implicit reductions with branch and conditioning constraints. These reductions are specialized to power-law, Hill-type, and polynomial growth laws to identify when sparse regression can recover mechanistic parameters and when low-residual fits are nonidentifiable or nontransferable.
\end{abstract}



\section{Introduction}
Many biological populations contain subpopulations that differ in proliferative status. At the single-cell or subpopulation level, cells can actively cycle and divide, withdraw into a reversible quiescent or dormant state, and later re-enter proliferation when environmental, developmental, or therapeutic conditions change \cite{yao2014quiescence,li2010quiescent}. This behavior appears across several biological systems. Active and quiescent adult stem-cell populations can coexist within tissues \cite{li2010quiescent}. Naive T cells are maintained in quiescence and exit quiescence after antigen and costimulatory signals before clonal expansion \cite{chapman2018hallmarks}. Bacterial persister cells can switch between normal growth and a slow-growing drug-tolerant phenotype \cite{balaban2004bacterial,kussell2005bacterial}. Quiescent cancer cells can evade therapies that primarily target cycling cells and later re-enter proliferation to contribute to recurrence \cite{lindell2023quiescent}. These biological observations motivate active-quiescent compartment models, by which we mean models that partition a population into a proliferative active compartment and a nonproliferative quiescent compartment, with transitions into and out of quiescence. Variants of this active-quiescent structure have been used to model quiescent phases in mathematical biology, microbial persistence, stem-cell renewal, immune activation, and tumor response to antimitotic therapy \cite{hadeler2008quiescent,kussell2005bacterial,kozusko2007microenvironment}. In these models, growth acts primarily through the active compartment, while the quiescent compartment stores cells that may later return to the active state and may be hidden or only indirectly observed. This distinction is biologically important because changes in switching can alter the apparent growth of the measured population without changing the intrinsic proliferation law of active cells. For example, chemotherapy may reduce proliferating tumor cells while sparing quiescent cells, changing the observed tumor trajectory. Conversely, reactivation of quiescent stem, immune, microbial, or cancer subpopulations can produce apparent population expansion even if the active-state growth law is unchanged. Thus, observed population dynamics reflect both proliferation within the active state and movement between active and quiescent states on experimentally relevant time scales.

The measurements available in applications rarely reveal every compartment of the model. In cell-population experiments, active-state information may come from proliferation or cell-cycle markers such as Ki-67, BrdU/EdU incorporation, FUCCI reporters, phospho-histone H3 staining, or other marker panels measured by microscopy or flow cytometry \cite{scholzen2000ki67,salic2008edu,sakaue2008fucci}. These measurements may be sparse in time when marker panels are limited, when cells must be fixed or destructively processed, or when active-state assays are combined with separate population-level readouts. Other experiments provide primarily aggregate measurements, such as total nuclei counts, confluence, optical density, viable biomass, aggregate fluorescence, or tumor volume. Such readouts are common in microbial growth and persistence studies, cell-culture assays, and tumor growth experiments, including models with quiescent or drug-tolerant compartments \cite{balaban2004bacterial,kussell2005bacterial,kozusko2007microenvironment,brown2017tgf}. These measurements combine active and quiescent states and therefore do not directly reveal the compartment in which growth occurs. Thus the observable determines which differential equation is available to the data after hidden variables have been eliminated.

The focus of this work is to analyze how active-quiescent structure and partial biological observation affect regression-based equation learning. Sparse regression methods such as SINDy and PDE-FIND infer parsimonious differential models by selecting a small number of active terms from a prescribed candidate library \cite{brunton2016sindy,rudy2017pdefind}. Recent work has addressed partial observation using higher-derivative dictionaries, encoder-based hidden-state reconstruction, delay-coordinate autoencoders, and sparse data assimilation \cite{somacal2022hidden,lu2022partial,bakarji2023partial,ribera2022model}. These methods show that hidden variables are a central difficulty in equation discovery. Here we take a complementary, constructive viewpoint: observable reduction can be used to derive the equation satisfied by the measured signal and to translate that equation into a sparse-regression target, candidate library, and coefficient diagnostics. Active-marker and aggregate readouts enter because they generally induce different observable reductions, so the same hidden active-quiescent mechanism may require different sparse-regression libraries. The goal is therefore to use observable reduction to guide sparse regression toward mechanistically interpretable fits instead of relying on a generic library in the measured variables.

This dependence creates a structural difficulty in hidden compartment biological systems. In the standard SINDy formulation, time series measurements are assembled into a data matrix \(X\), derivative data are assembled into \(\dot{X}\), and one seeks a sparse coefficient matrix \(\Xi\) satisfying
\[
\dot{X}=\Theta(X)\Xi,
\]
where \(\Theta(X)\) is a library of candidate nonlinear functions evaluated on the measured state. If the full active-quiescent state $x$ were observed, this framework would suggest constructing \(X\) from measurements of the full state, and fitting a sparse model for the vector field in the full state space. In many biological experiments, however, not all state variables are observable. The available data may contain only an active-state marker, an aggregate population measurement, or sparse marker measurements combined with aggregate readouts. In such cases, SINDy cannot be applied directly to the full active-quiescent state. A library built only from the observed variables may still fit the measured trajectory, but it might not recover the hidden switching structure, the active growth law, or a mechanism that transfers across experimental regimes.

We refer to the differential equation obtained after eliminating hidden states as the observable reduction. In the two compartment setting studied here, this reduction gives an equation for the observable state \(y\) of the form
\[
\ddot{y}=g(y,\dot{y}).
\]
The function \(g\) depends on which state is observed. For an active-state observable, the hidden quiescent variable can often be eliminated to obtain a closed equation involving only the active state and its derivatives \cite{hadeler2008quiescent}. In this case, the observable reduction can preserve a direct relationship between discovered coefficients and the switching and growth parameters. For an aggregate observable, the measured state combines active and quiescent populations. The reduction then requires reconstructing the hidden active state from the aggregate growth relation. Consequently, two datasets generated by the same hidden active and quiescent model can require different discovery libraries because they correspond to different observed states.

In this work, we show that sparse regression for hidden-state systems should begin with the equation satisfied by the measured quantity. For active-quiescent models, this means deriving the scalar observable reduction before choosing a candidate library. A flexible library may reproduce one observed trajectory with small residual error while producing coefficients that do not map to switching or growth parameters, violate the active-quiescent reduction, or fail outside the training regime. We therefore use observable reduction as a structural check on both the selected terms and the fitted coefficients. We investigate this idea for nonlinear active-quiescent ODE models with active growth and linear switching between active and quiescent states. We consider two observation maps: direct observation of the active state and observation of the aggregate population. For each map, we derive the corresponding scalar observable reduction and specialize it to power-law, Hill-type, and polynomial active growth. These growth laws produce different sparse-regression and identifiability issues. Power-law growth gives structured libraries with explicit coefficient-to-parameter maps. Hill-type growth introduces shape parameters whose practical identifiability depends on the observed dynamical regime. Polynomial growth produces coefficient-pairing constraints under active observation and inverse branch ambiguity under aggregate observation.

Our contributions are as follows. First, we derive scalar observable reductions for active-state and aggregate observations of nonlinear active-quiescent ODE models. The active-state reduction builds on existing active-quiescent reduction results \cite{hadeler2008quiescent}. The aggregate-state reduction identifies the inverse reconstruction, branch-selection, and conditioning issues introduced by population-level measurements. Second, we apply the observable reductions to common biological growth laws, including power-law, Hill-type, and polynomial active growth. This gives observation-specific libraries, coefficient maps, and consistency diagnostics. Third, we use synthetic computational benchmarks to compare the reduction-guided libraries with polynomial SINDy surrogates and to identify failures caused by partial or aggregate observation. These failures include low-residual but noninterpretable fits, transfer failure, branch ambiguity, and regime-dependent nonidentifiability.

\section{Active-Quiescent Model and Observation Maps}\label{sec:model}

In this section, we describe a two-compartment active-quiescent model with two scalar observation maps: active-state observation and aggregate population observation. The same hidden dynamics induce different scalar equations depending on which quantity is measured.

\subsection{Hidden active-quiescent dynamics}\label{sec:hidden-model}

We consider the active and quiescent ODE model \cite{hadeler2008quiescent}
\begin{subequations}\label{eq:model}
\begin{align}
    \dot v &= f(v)-\gamtwo v+\gamone w, \label{eq:model-v}\\
    \dot w &= -\gamone w+\gamtwo v. \label{eq:model-w}
\end{align}
\end{subequations}
Here \(v(t)\) denotes the active population, \(w(t)\) denotes the quiescent population, and \(f(v)\) is the intrinsic active growth law. The parameter \(\gamtwo\geq 0\) is the active to quiescent switching rate, and \(\gamone\geq 0\) is the quiescent to active switching rate. Growth acts only through the active compartment, while switching redistributes mass between the two compartments. Thus net population growth is controlled by the active state through \(f(v)\). This identity is central for aggregate observations because \(u=v+w\) satisfies
\begin{equation}\label{eq:u-dot-prelim}
    \dot u=f(v).
\end{equation}
An aggregate measurement therefore depends on the hidden active state even when only the total population is observed.

The two switching rates have distinct dynamical roles. Increasing \(\gamtwo\) transfers active mass into the quiescent compartment and can slow the apparent growth of the active population. Increasing \(\gamone\) transfers quiescent mass back into the active compartment and can restore growth by feeding the proliferative state. When both rates are positive, switching is reversible. When one rate is zero, the model reduces to a one-way transition.

We do not impose a specific form for \(f(v)\) at this stage. The reductions below first treat \(f\) as a general nonlinear active growth law and then specialize to {power-law}, {Hill-type}, and polynomial growth in Section~\ref{sec:libraries}. This separation is important because the hidden model specifies where growth occurs, while the analytic form of \(f\) determines the nonlinear terms that appear in the scalar observable equation.

\subsection{Observation maps}\label{sec:observation-maps}

The hidden model \eqref{eq:model} is first order in the two-dimensional state
\(x(t)=(v(t),w(t))^T\), but the measured variable is not always the full state.
We represent a scalar measurement by an observation map
\begin{equation}\label{eq:general-observation-map}
    y(t)=P(x(t))=P(v(t),w(t)).
\end{equation}
This notation allows for general nonlinear measurement processes. In the main
analysis, we focus on two linear observation maps of the form
\begin{equation}\label{eq:observation-maps}
    y(t)=C x(t),
    \qquad
    C \in \mathbb{R}^{1\times 2}.
\end{equation}
Active-state observation corresponds to
\[
    C_{\mathrm{act}}=(1,0),
    \qquad
    y(t)=C_{\mathrm{act}}x(t)=v(t),
\]
while aggregate population observation corresponds to
\[
    C_{\mathrm{agg}}=(1,1),
    \qquad
    y(t)=C_{\mathrm{agg}}x(t)=v(t)+w(t)=u(t).
\]
The first linear map represents active-state measurements, such as active markers, cell-cycle indicators, proliferating cell counts, or reporter signals associated with active growth. The second represents aggregate measurements, such as total population size, cell density, optical density, fluorescence intensity, viable biomass, or tumor burden. Under active observation, the measured variable is one component of the hidden state, so the quiescent variable can be eliminated from the active equation. Under aggregate observation, the measured variable combines active and quiescent states, so the active state must be reconstructed from $\dot u=f(v)$. If \(f\) has a locally admissible inverse branch \(F\), then
\begin{equation}\label{eq:inverse-growth}
    v=F(\dot u),
    \qquad
    f(F(\dot u))=\dot u.
\end{equation}
For nonlinear \(f\), this inverse relation may be local, multivalued, or valid only in restricted dynamical regimes.

\begin{figure}[t]
\centering
\includegraphics[width=0.90\linewidth]{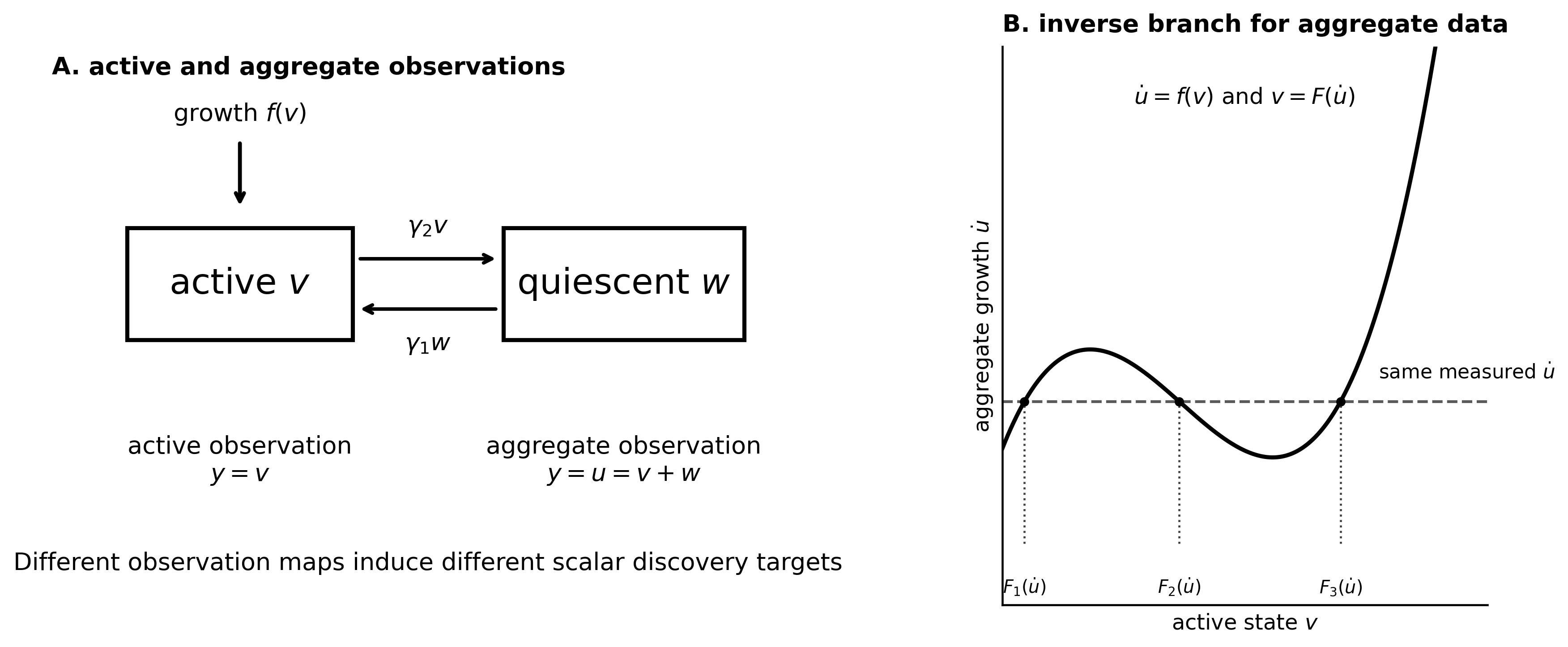}
\caption{Conceptual active-quiescent model, observation maps, and aggregate inverse branches. A) Active growth \(f(v)\) occurs only in the active compartment, with switching between active and quiescent states. Active observation measures \(y=v\) and leads to a closed observable reduction for the measured active state. Aggregate observation measures \(y=u=v+w\), so the hidden active state must be reconstructed from the aggregate growth relation \(\dot u=f(v)\). B) The same measured value of \(\dot u\) can intersect a nonmonotone growth curve at more than one positive value of \(v\). In that case, the inverse reconstruction is branch dependent, and aggregate data alone do not determine which active state generated the observed growth.}
\label{fig:schematic}
\end{figure}

Figure~\ref{fig:schematic} summarizes the hidden model and the two observation maps. The corresponding scalar sparse-regression targets are stated in Section~\ref{sec:reductions}, with full derivations given in Appendix~\ref{app:observable-derivations}.

\subsection{Synthetic cell-population measurement example}

Consider a computationally generated cell-population experiment in which cells occupy either an active proliferative state or a quiescent nonproliferative state. The hidden variables are the active population \(v(t)\) and the quiescent population \(w(t)\). Active cells grow according to \(f(v)\), switch to quiescence at rate \(\gamma_2\), and return from quiescence at rate \(\gamma_1\). An active-state marker, such as a proliferation reporter, cell-cycle marker, or active-state fluorescence signal, gives a measurement close to \(y(t)=v(t)\). A total-population readout, such as nuclei count, confluence, optical density, viable biomass, tumor burden, or aggregate fluorescence, gives a measurement close to \(y(t)=u(t)=v(t)+w(t)\). Although both measurements come from the same hidden active-quiescent mechanism, they lead to different scalar regression targets. Active-state observation gives the second-order active observable reduction. Aggregate observation instead requires reconstructing the hidden active state from \(\dot u=f(v)\). The measurement process therefore determines the regression target before sparse regression is applied. Table~\ref{tab:bio-observation-maps} summarizes the corresponding observables, regression targets, and diagnostics.

\begin{table}[tbp]
\centering
\small
\begin{tabular}{p{0.25\textwidth}p{0.16\textwidth}p{0.25\textwidth}p{0.24\textwidth}}
\toprule
Biological measurement & Observable & Regression target & Primary diagnostic \\
\midrule
Active marker or proliferation reporter & \(y=v\) & second-order active reduction & coefficient maps and coefficient pairing \\
Total nuclei, confluence, OD, biomass, tumor burden, or aggregate fluorescence & \(y=u=v+w\) & aggregate inverse or implicit reduction & positive branch count and inverse conditioning \\
Sparse active-marker measurements plus aggregate readout & \(u(t), v(t_i)\) & aggregate reduction with active-state constraints & branch resolution and parameter constraints \\
Perturbed aggregate readout & \(u_{\mathrm{pert}}(t)\) & shared mechanism across conditions & transfer or perturbation prediction error \\
\bottomrule
\end{tabular}
\caption{Biological interpretation of observation maps. Different biological measurement protocols correspond to different mathematical observables, sparse-regression targets, and identifiability diagnostics.}
\label{tab:bio-observation-maps}
\end{table}

\section{Sparse Identification of Nonlinear Dynamics}\label{sec:sindy-library-problem}

Sparse identification of nonlinear dynamics (SINDy) assumes that an observed vector field can be represented as a sparse expansion in prescribed candidate functions \cite{brunton2016sindy,kaptanoglu2022pysindy,messenger2021weak}. Let \(\{\theta_j\}_{j=1}^p\) be a library of candidate nonlinear functions. For a measured trajectory \(x(t)\in\mathbb{R}^d\), SINDy seeks sparse coefficients \(\xi_{ji}\) such that
\[
    \dot{x}_i(t) \approx \sum_{j=1}^p \xi_{ji}\theta_j(x(t)),
    \qquad i=1,\ldots,d.
\]
At sampled times \(t_m\), this gives the matrix regression problem
\[
    \dot{X}\approx \Theta(X)\Xi,
\]
where \(\Theta(X)_{mj}=\theta_j(x(t_m))\) and \(\Xi=(\xi_{ji})\). More generally, the sparse regression step does not require the left hand side to be a first derivative of the full state. One may choose the regression target from the differential equation available for the measured quantity. In standard first-order SINDy, this target is \(\dot{X}\). A major issue in biological applications is that not all state variables are observable. In many experiments, the available data contain only a subset of states or aggregate quantities. The sparse regression problem is therefore coupled to the observation process. A library built only from the observed variables may fit the data well while missing hidden switching dynamics, the active growth law, or a mechanism that transfers across regimes. In addition, choosing a library with correlated functions can lead to ill-conditioning and model misidentification, especially when biological trajectories sample a restricted region of state space \cite{feng2026illconditioning}. In an observable-reduction setting, hidden states are first eliminated, and the resulting scalar equation may instead use the second derivative of the observable as the regression target. Thus a reduced equation of the form \(\ddot{y}=g(y,\dot{y})\) leads to a sparse regression problem in which samples of \(\ddot{y}\) are fit using candidate functions built from \(y\) and \(\dot{y}\).

\section{Observable Reductions}\label{sec:reductions}

The observable reduction is the second-order derivative equation satisfied by the measured variable after hidden states have been eliminated. These reductions provide the regression targets and guide the library construction used in the following sections. The derivations are described in Appendix~\ref{app:observable-derivations}. 

\subsection{Active-state observable reduction}

Assume \(f\in C^1\) and \(\gamone>0\). If \(y=v\), then the active state satisfies the closed second-order equation
\begin{equation}\label{eq:active-general}
    \boxed{
    \vddot=\bigl(f'(v)-\gamone-\gamtwo\bigr)\vdot+\gamone f(v).
    }
\end{equation}
Thus active observation does not lead to a first-order scalar equation in general. The hidden quiescent compartment appears through the acceleration term and through the switching parameters \(\gamone\) and \(\gamtwo\). For sparse regression, the appropriate target is therefore a second-order model \cite{hadeler2008quiescent} of the form
\begin{equation}\label{eq:sindy-active-target}
    \ddot v=\ThetaLib(v,\dot v)\xi.
\end{equation}
The candidate library should contain the terms generated by \(f(v)\), by \(f'(v)\dot v\), and by the linear damping term \(\dot v\).

The assumption \(\gamone>0\) is used when eliminating \(w\) from the active equation. If \(\gamone=0\), quiescent cells do not return to the active compartment and the active equation is already first-order closed:
\[
    \dot v=f(v)-\gamtwo v.
\]
In that one-way switching case, a first-order active-state regression target may be appropriate. Thus the regression target depends on both the observation map and whether switching is reversible.

The active reduction is relatively favorable for mechanistic interpretation. When the correct growth family is used, the coefficients of \eqref{eq:active-general} can often be mapped back to switching rates and growth parameters. This is the basis for the coefficient-to-parameter maps and consistency diagnostics developed in Section~\ref{sec:libraries}.

\subsection{Aggregate-state observable reduction}

If \(y=u\), then the aggregate variable satisfies (\ref{eq:u-dot-prelim}). The measured aggregate growth rate is therefore the active growth law evaluated at the hidden active state. The corresponding differential algebraic reduction is
\begin{align} \label{eq:agg-dae}
    \uddot &= f'(v)\left[\udot+\gamone u-(\gamone+\gamtwo)v\right].
\end{align}
If \(f\) has a locally admissible inverse branch \(F\), then \(v=F(\udot)\) and the aggregate variable satisfies the explicit scalar equation
\begin{equation}\label{eq:agg-general}
    \boxed{
    \uddot=f'\!\left(F(\udot)\right)
    \left[\udot+\gamone u-(\gamone+\gamtwo)F(\udot)\right].
    }
\end{equation}
Unlike the active reduction, the aggregate reduction depends on the inverse structure of \(f\). If \(f\) is not one-to-one on the biological domain, then different positive active states may produce the same value of \(\dot u\). The aggregate equation is then branch dependent. This branch dependence is illustrated schematically in Figure~\ref{fig:schematic}B, where a single measured value of \(\dot u\) can intersect the growth curve at one or multiple positive active states.

For sparse regression, aggregate observation is therefore more restrictive than active observation. A standard library in \((u,\dot u)\) is appropriate only when the inverse branch can be represented explicitly and stably. Otherwise, the discovery problem is better treated with an implicit or nonlinear formulation based on \eqref{eq:agg-dae}. This distinction is central for the aggregate Hill and aggregate polynomial cases below.

\section{Observable-Reduction-Guided Sparse-Regression Libraries}\label{sec:libraries}

We focus on the observation-specific libraries used in the numerical experiments. Substituting the power-law, Hill-type, and polynomial growth laws into the active and aggregate reductions gives the corresponding regression targets, candidate libraries, coefficient maps, and structural diagnostics. The detailed algebraic substitutions are described in Appendix~\ref{app:growth-library-derivations}.

\mbox{}
\subsection{Power-law growth}

Let
\begin{equation}\label{eq:power-f}
    f(v)=av^n,\qquad a>0,\qquad n>0.
\end{equation}
For active observation, the second-order equation is
\begin{equation}\label{eq:active-power}
    \boxed{
    \vddot=anv^{n-1}\vdot-(\gamone+\gamtwo)\vdot+\gamone av^n.
    }
\end{equation}
For known \(n\), the active observable library is
\begin{equation}\label{eq:active-power-library}
    \mathcal{L}^{\mathrm{act}}_{\mathrm{power}}
    =
    \left\{v^{n-1}\vdot,\;\vdot,\;v^n\right\}.
\end{equation}
If
\[
    \vddot=c_1v^{n-1}\vdot+c_2\vdot+c_3v^n,
\]
then the mechanistic parameters are recovered from
\begin{equation}\label{eq:active-power-map}
    a=\frac{c_1}{n},\qquad
    \gamone=\frac{nc_3}{c_1},\qquad
    \gamtwo=-c_2-\gamone.
\end{equation}
This coefficient-to-parameter map is conditional on the assumed power-law growth family and on the exponent \(n\) being known or correctly selected. It also requires that the corresponding library terms are present and that the coefficient \(c_1=an\) is nonzero and reliably estimated. If the selected exponent is incorrect, a sparse fit may still have a small residual on a restricted trajectory, but the resulting coefficients no longer have the stated mechanistic interpretation.

For aggregate observation, the second-order equation is
\begin{equation}\label{eq:agg-power}
    \boxed{
    \uddot
    =n a^{1/n}(\udot)^{(2n-1)/n}
    +n\gamone a^{1/n}u(\udot)^{(n-1)/n}
    -n(\gamone+\gamtwo)\udot.
    }
\end{equation}
The aggregate power-law reduction uses the positive inverse branch \(v=(\dot u/a)^{1/n}\). For non-integer \(n\), this real-valued branch requires \(\dot u>0\), corresponding to positive active growth on the observed trajectory. Therefore the aggregate power-law library is valid on trajectory regions where the observed aggregate derivative remains positive and the positive active-state branch is biologically admissible.
For known \(n\), the aggregate observable library is
\begin{equation}\label{eq:agg-power-library}
    \mathcal{L}^{\mathrm{agg}}_{\mathrm{power}}
    =
    \left\{(\udot)^{(2n-1)/n},\;u(\udot)^{(n-1)/n},\;\udot\right\}.
\end{equation}
If
\[
    \uddot=c_1(\udot)^{(2n-1)/n}+c_2u(\udot)^{(n-1)/n}+c_3\udot,
\]
then
\begin{equation}\label{eq:agg-power-map}
    a=\left(\frac{c_1}{n}\right)^n,\qquad
    \gamone=\frac{c_2}{c_1},\qquad
    \gamtwo=-\frac{c_3}{n}-\gamone.
\end{equation}

\subsection{Hill-type growth}

Let
\begin{equation}\label{eq:hill-f}
    f(v)=\frac{av^n}{K^n+v^n},
    \qquad a>0,\quad K>0,\quad n>0.
\end{equation}
For active observation, the second-order equation is
\begin{equation}\label{eq:active-hill}
    \boxed{
    \vddot=
    \left[an\frac{K^n v^{n-1}}{(K^n+v^n)^2}-(\gamone+\gamtwo)\right]\vdot
    +\gamone a \frac{v^n}{K^n+v^n}.
    }
\end{equation}
For known \((n,K)\), the active Hill library contains the three functions
\begin{equation}\label{eq:active-hill-atoms}
    D_{n,K}(v,\vdot)=\frac{K^n v^{n-1}}{(K^n+v^n)^2}\vdot,
    \qquad
    G_{n,K}(v)=\frac{v^n}{K^n+v^n},
    \qquad
    \vdot.
\end{equation}
The corresponding coefficients satisfy \(c_D=an\), \(c_G=\gamone a\), and \(c_{\dot v}=-(\gamone+\gamtwo)\). Thus the active Hill problem is linear in the coefficients once the shape parameters \((n,K)\) are fixed. 

For aggregate observation, the inverse relation is valid on the positive biological domain when
\begin{equation}\label{eq:hill-q-condition}
    0<\dot u<a.
\end{equation}
On this admissible range, the hidden active state can be written as
\begin{equation}\label{eq:hill-inverse}
    v=K\left(\frac{\dot u}{a-\dot u}\right)^{1/n}.
\end{equation}
The aggregate Hill target can then be written as
\begin{equation}\label{eq:agg-hill}
    \boxed{
    \uddot=
    \frac{n\dot u(a-\dot u)}{av}
    \left[\dot u+\gamone u-(\gamone+\gamtwo)v\right],
    \qquad
    v=K\left(\frac{\dot u}{a-\dot u}\right)^{1/n}.
    }
\end{equation}
This is not an ordinary linear SINDy problem in a fixed library because the inverse reconstruction depends on unknown parameters. The aggregate Hill case therefore requires nonlinear search or an implicit formulation together with the biological constraint \(0<\dot u<a\). In particular, admissibility is parameter dependent. A candidate value of \(a\) is valid only if \(0<\dot u(t)<a\) for all time points used in the fit. Candidate parameter sets violating this inequality do not define a positive real active-state reconstruction and should be excluded or penalized.

The Hill law also gives practical identifiability regimes. If \(v^n\ll K^n\), then \(f(v)\approx(a/K^n)v^n\), so data mainly identify the lumped ratio \(a/K^n\). If \(v^n\gg K^n\), then \(f(v)\approx a\) and \(f'(v)\approx0\), so \(K\) has little influence on the observed dynamics. Reliable separation of \(a\), \(K\), and \(n\) requires trajectories that sample the transition region \(v\approx K\).

\subsection{Polynomial growth}

Let
\begin{equation}\label{eq:poly-f}
    f(v)=\sum_{j=m}^{p}a_jv^j,
    \qquad 1\leq m\leq p,
    \qquad a_p\neq 0.
\end{equation}
For active observation, the second-order equation is
\begin{equation}\label{eq:active-poly}
    \boxed{
    \vddot=
    \sum_{j=m}^{p}j a_jv^{j-1}\vdot
    -(\gamone+\gamtwo)\vdot
    +\gamone\sum_{j=m}^{p}a_jv^j.
    }
\end{equation}
The active polynomial library is
\begin{equation}\label{eq:active-poly-library-correct}
    \mathcal{L}^{\mathrm{act}}_{\mathrm{poly}}
    =
    \left\{v^{m-1}\vdot,\ldots,v^{p-1}\vdot,\;\vdot,\;v^m,\ldots,v^p\right\}.
\end{equation}
The key structural feature is coefficient pairing. The coefficient of \(v^{j-1}\vdot\) is \(j a_j\), while the coefficient of \(v^j\) is \(\gamone a_j\). Therefore, whenever \(a_j\neq0\) and both coefficients are estimated,
\begin{equation}\label{eq:gamma1-pair}
    \gamone
    =
    j\frac{\text{coefficient of }v^j}{\text{coefficient of }v^{j-1}\dot v}.
\end{equation}
Consistency of \eqref{eq:gamma1-pair} across multiple degrees is a structural signature of hidden active and quiescent dynamics. It tests whether the fitted polynomial acceleration obeys the coefficient relations imposed by eliminating the hidden quiescent state.

For aggregate observation, the polynomial reduction is most naturally retained in differential algebraic form
\begin{subequations}\label{eq:agg-poly-dae}
\begin{align}
    \udot &=\sum_{j=m}^{p}a_jv^j,\\
    \uddot
    &=
    \left(\sum_{j=m}^{p}j a_jv^{j-1}\right)
    \left[\udot+\gamone u-(\gamone+\gamtwo)v\right].
\end{align}
\end{subequations}
Recovering \(v\) requires solving
\begin{equation}\label{eq:poly-root}
    \sum_{j=m}^{p}a_jv^j-\udot=0.
\end{equation}
If this equation has a unique positive real solution over the observed range, then the hidden active state can be reconstructed on that biological branch, assuming the polynomial growth law and its coefficients are known. If two or more positive roots exist for the same value of \(\dot u\), aggregate data alone cannot determine which active state generated the observation. Near \(f'(v)=0\), the inverse branch becomes ill-conditioned because \(d v/d\dot u=1/f'(v)\). This branch-reconstruction condition should not be interpreted as full structural identifiability of all unknown coefficients in \(f\) and switching rates \(\gamone,\gamtwo\). Parameter identifiability additionally requires enough independent information in the reconstructed branch and aggregate equation to determine those unknowns.

\subsection{Observable-reduction-induced general library construction for active observable}
The power-law, Hill, and polynomial libraries above are controlled benchmarks that test recovery when a growth family is specified. More generally, the active observable reduction supplies a construction rule for a dictionary of candidate growth atoms: each candidate $\phi(v)$ contributes the pair $\phi(v)$ and $\phi'(v)\dot v$, together with the shared damping column $\dot v$. We therefore use the duplicate-free 26-term library
\begin{equation}\label{eq:general-active-library}
\begin{aligned}
\mathcal L_{\mathrm{act}}={}&\{\dot v\}
\cup\{v,v^2,v^3,v^4\}
\cup\{2v\dot v,3v^2\dot v,4v^3\dot v\}\\
&\cup\left\{G_{n,K}(v),G'_{n,K}(v)\dot v:
 n\in\{1,2,3\},\ K\in\{0.5,1,2\}\right\},
\end{aligned}
\end{equation}
where
\begin{equation}\label{eq:general-active-hill-atoms}
G_{n,K}(v)=\frac{v^n}{K^n+v^n},\qquad
G'_{n,K}(v)\dot v=\frac{nK^n v^{n-1}}{(K^n+v^n)^2}\dot v.
\end{equation}

Integer power-law active libraries are nested in the polynomial block because an atom $v^n$ induces both $v^n$ and $n v^{n-1}\dot v$. In particular, the controlled power-law case with $n=2$ is represented by $v^2$ and $2v\dot v$ together with the shared $\dot v$ term. The controlled polynomial case with maximum degree $p=3$ is also contained in the candidate polynomial block with $p=4$. Therefore the 26-term dictionary contains the controlled power-law, polynomial, and Hill active-observable libraries. 

Table \ref{tab:libraries} provides a summary for the observable-reduction libraries and coefficients relations.

\begin{table}[t]
\centering
\small
\caption{Observable-reduction libraries and coefficient relations.}
\label{tab:libraries}
\begin{tabular}{p{0.13\linewidth}p{0.13\linewidth}p{0.3\linewidth}p{0.3\linewidth}}
\toprule
Growth law & Observation & Library terms & Identifiable relations \\
\midrule
Power \(av^n\) & \(y=v\) &
\(v^{n-1}\dot v,\dot v,v^n\) &
\(a=c_1/n\), \(\gamma_1=nc_3/c_1\), \(\gamma_2=-c_2-\gamma_1\) \\
Power \(av^n\) & \(y=u\) &
\(\dot u^{(2n-1)/n}, u\dot u^{(n-1)/n},\dot u\) &
\(a=(c_1/n)^n\), \(\gamma_1=c_2/c_1\), \(\gamma_2=-c_3/n-\gamma_1\) \\
\hline
Hill & \(y=v\) &
\(D_{n,K}(v,\dot v), \dot v, G_{n,K}(v)\) &
Linear coefficients after fixing \(n,K\) \\
Hill & \(y=u\) &
Implicit or inverse terms involving \(v=K(\dot u/(a-\dot u))^{1/n}\) &
Nonlinear parameter search. Valid for \(0<\dot u<a\) \\
\hline
Polynomial & \(y=v\) &
\(v^{j-1}\dot v,\dot v,v^j\) &
\(\gamma_1=j\,c(v^j)/c(v^{j-1}\dot v)\) across degrees \\
Polynomial & \(y=u\) &
Positive branch of \(f(v)=\dot u\) &
Unique positive branch required for aggregate identifiability \\
\bottomrule
\end{tabular}
\end{table}

\section{Numerical Methods}\label{sec:numerics}

We generated synthetic trajectories from \eqref{eq:model} and fitted candidate sparse-regression models using the targets and libraries derived above. We then evaluated each fit using the regression targets, scoring metrics, and validation checks described below.

\subsection{Regression strategies}\label{sec:discovery-strategies}

The regression strategy depends on the observation map. For active observations, we compare three approaches. The first fits a naive first-order polynomial model for \(\dot v\), even though reversible active-quiescent switching generally does not produce a closed first-order equation for \(v\). The second fits \(\ddot v\) with a generic polynomial library in \((v,\dot v)\). The third fits \(\ddot v\) with the 26-term library induced by the active observable reduction. This dictionary contains polynomial atoms and candidate Hill atoms together with their derivative-coupled partners. It therefore tests whether the growth structure can be selected without assuming a known growth family.

These approaches answer different questions. The first-order fit tests whether a trajectory-dependent surrogate can reproduce \(\dot v\). The generic second-order fit tests whether a flexible polynomial basis can reproduce the correct acceleration target. The 26-term fit additionally tests support recovery and permits the selected coefficients to be interpreted using the growth-specific maps in Section~\ref{sec:libraries}. For aggregate observations, the inverse dependence of the observable reduction prevents one common fixed dictionary from playing the same role. We therefore retain growth-specific inverse or nonlinear procedures for the aggregate power-law and aggregate Hill cases.

Residual error is reported for all fits, but it is not treated as evidence of mechanism recovery by itself. A low-residual model may target the wrong derivative, select terms without a valid coefficient map, violate active-quiescent coefficient relations, or fail when evaluated in another dynamical regime. We therefore interpret residuals together with selected support, recovered parameters, coefficient-pairing consistency, and transfer error.

\subsection{Synthetic data generation and implementation details}

Synthetic trajectories were generated from \eqref{eq:model}. The simulator stores \(v,w,u\) and the derivatives \(\dot v,\dot u,\ddot v,\ddot u\) computed from the model and observable reductions. The present experiments use exact synthetic trajectories and model-evaluated derivatives by design. Our goal is to isolate structural failures caused by partial and aggregate observation from errors caused by finite sampling, derivative estimation, or measurement noise. These experiments therefore provide a best-case identifiability benchmark for the observation protocol. Robustness to noise is a separate practical-identifiability question that we leave to future extensions using weak and implicit formulations.
Matched active and aggregate controlled benchmarks use the same underlying growth choices: the power-law exponent is \(n=2\), the polynomial maximum degree is \(p=3\), and the Hill truth is \(a=1.5\), \(K=1\), and \(n=2\). 


\subsection{Scoring, validation, and reproducibility}

For a chosen library \(\ThetaLib\), the regression problem is
\begin{equation}\label{eq:stlsq}
    \min_{\xi}\|\mathbf z-\ThetaLib\xi\|_2^2+\lambda\|\xi\|_2^2
    \quad\text{followed by coefficient thresholding,}
\end{equation}
where \(\mathbf z\) is the target derivative. For the first-order derivative equations, \(\mathbf z\) is either \(\dot v\) or \(\dot u\). For the second-order derivative equations, \(\mathbf z\) is either \(\ddot v\) or \(\ddot u\). Before sparse regression, each library column is normalized to unit Euclidean norm. If \(\ThetaLib=[\Theta_1,\ldots,\Theta_p]\), we fit with scaled columns \(\widetilde \Theta_j=\Theta_j/\|\Theta_j\|_2\), with zero-norm columns left unchanged. Thresholding is performed in physical coefficient units, and the final coefficients are transformed back by \(\widehat\xi_j=\widetilde\xi_j/\|\Theta_j\|_2\). All coefficient-to-parameter maps are therefore applied to unscaled coefficients. The reported residual is
\begin{equation}\label{eq:rel-resid}
    \mathcal R=\frac{\|\mathbf z-\ThetaLib\xi\|_2}{\|\mathbf z\|_2}.
\end{equation}
Residual error measures agreement with the observable regression target, but it does not by itself establish mechanistic recovery. For Hill-type growth, transfer across regimes is evaluated with relative target or prediction errors, for example \(\|z^{(s)}-\widehat z^{(s)}\|_2/\|z^{(s)}\|_2\) for a model trained in one regime and tested in regime \(s\).
For active polynomial recovery, we report the degree-wise consistency estimates
\[
    \widehat{\gamma}_{1,j}
    =
    j\frac{\widehat c(v^j)}{\widehat c(v^{j-1}\dot v)}.
\]
We summarize the spread of these degree-wise estimates using a normalized pairing-consistency score over the degrees for which both paired coefficients are selected and reliably nonzero. This score should be interpreted together with coefficient magnitudes and support recovery, since ratios involving very small denominator coefficients can be unstable.
For aggregate polynomial recovery, we diagnose the number of positive roots of \(f(v)-\dot u\) over the observed range and the conditioning factor \(|1/f'(v)|\). The root count diagnoses nonuniqueness of the aggregate inverse problem, while \(|1/f'(v)|\) diagnoses local numerical sensitivity of active-state reconstruction. A trajectory can therefore have a unique positive branch but still be practically unreliable if \(|f'(v)|\) is small along the observed range.

\section{Results}\label{sec:results}


\subsection{Active observations}
\subsubsection{Power-law recovery}

\begin{figure}[tbp]
\centering
\includegraphics[width=0.95\linewidth]{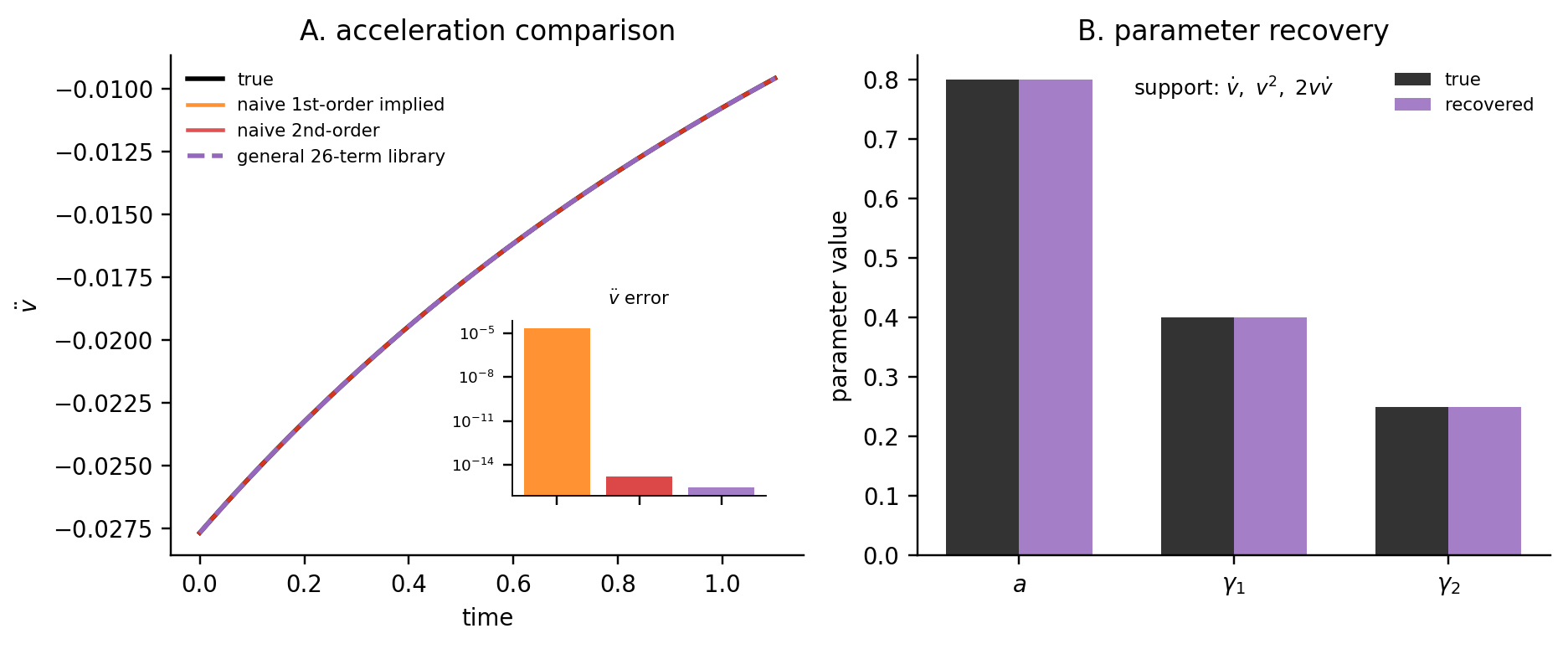}
\caption{Active power-law recovery from \(y=v\). Panel A compares the recovered accelerations and relative errors. Panel B reports parameter recovery and the selected 26-term-library support.}
\label{fig:active_power}
\end{figure}

For active-state observation with \(f(v)=0.8v^2\), the hidden quiescent state is eliminated only after differentiating the observed active state. The correct regression target is therefore the second-order equation \eqref{eq:active-power}, not a closed first-order law for \(v\). The naive first-order polynomial fit approximates \(\dot v\) along the sampled trajectory. Differentiating that surrogate gives the implied acceleration shown in Figure~\ref{fig:active_power}A, but its larger acceleration error confirms that fitting \(\dot v\) does not recover the active observable equation.

The generic second-order polynomial library and the 26-term library both contain the three terms required by the \(n=2\) observable reduction. Consequently, both recover the acceleration to numerical precision in this noise-free benchmark. Table~\ref{tab:active_power_equations} compares the true observable equation with the fitted equations. We found that the 26-term library selects the coherent support \(\{2v\dot v,\dot v,v^2\}\), which maps back to the correct switching parameters shown in Figure~\ref{fig:active_power}B.

\begin{table}[tbp]
\centering
\scriptsize
\input{tables/table_active_power_equations.tex}
\caption{Equation comparison for active power-law recovery.}
\label{tab:active_power_equations}
\end{table}

\subsubsection{Polynomial coefficient pairing}

The polynomial benchmark is more informative than the power-law case because the observable reduction predicts both the support and repeated coefficient relations. For \(f(v)=0.08v^2-0.01v^3\), the true acceleration contains the derivative-coupled terms \(0.16v\dot v\) and \(-0.03v^2\dot v\), the damping term \(-0.65\dot v\), and the growth terms \(0.032v^2\) and \(-0.004v^3\). The naive first-order polynomial again fits a different target. The full recovered first-order and second-order polynomial equations are reported in Table~\ref{tab:active-poly-equations}. The generic second-order library reproduces the acceleration equation because the true terms lie inside its polynomial basis.

\begin{figure}[tbp]
\centering
\includegraphics[width=0.95\linewidth]{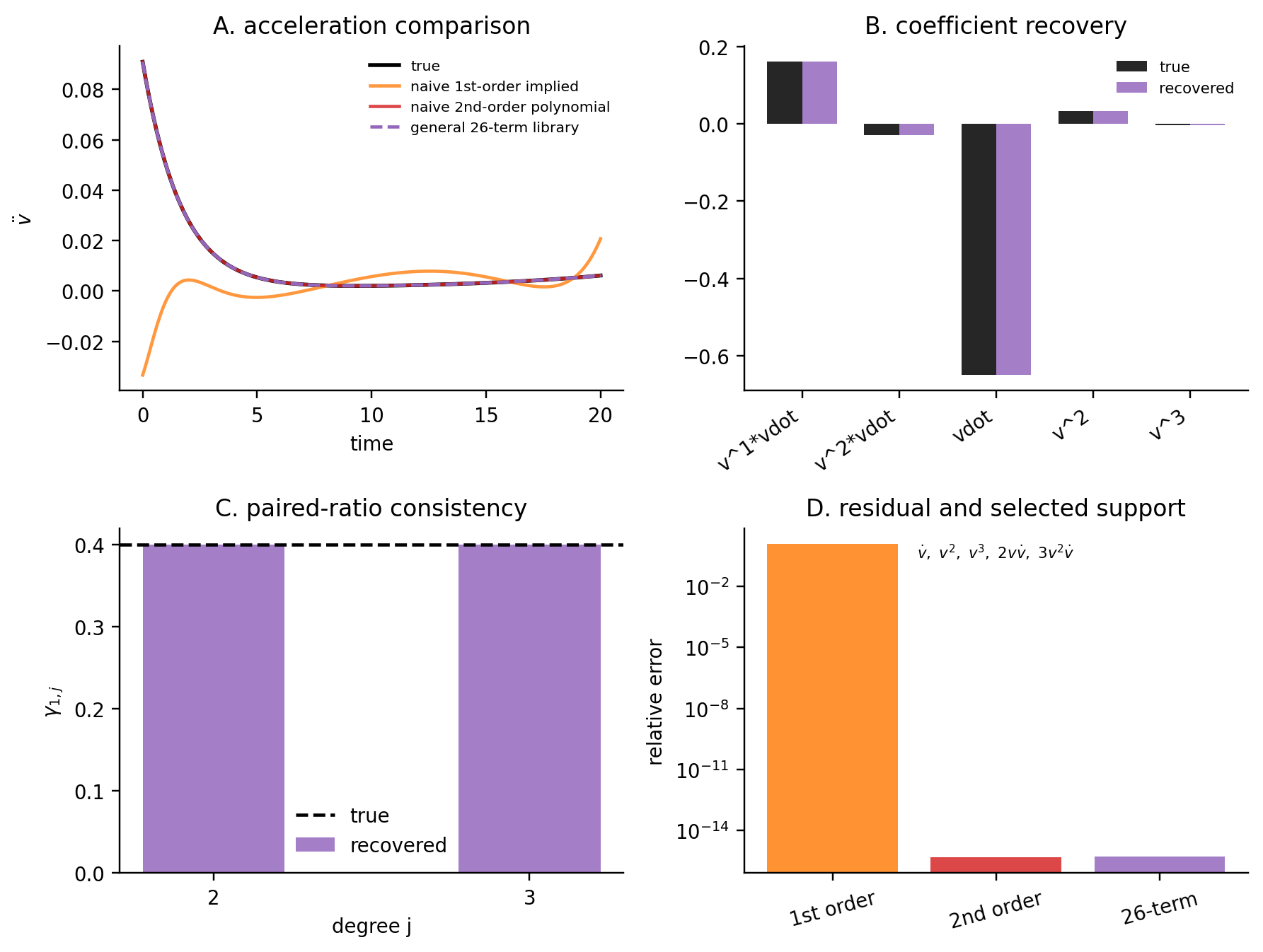}
\caption{Active polynomial recovery with the 26-term observable-reduction library. Panel A compares acceleration predictions. Panel B compares recovered and true coefficients. Panel C shows the degree-wise pairing estimates. Panel D reports residuals and selected support.}
\label{fig:active_poly_pairing}
\end{figure}

The 26-term library selects \(\dot v\), \(v^2\), \(v^3\), \(2v\dot v\), and \(3v^2\dot v\). Figure~\ref{fig:active_poly_pairing}B shows that the recovered observable-equation coefficients coincide with the true values. Their agreement in Figure~\ref{fig:active_poly_pairing}C is the mechanistic consistency check. Combining the common estimate \(\widehat\gamma_1=0.4\) with the coefficient of \(\dot v\) gives \(\widehat\gamma_2=0.65-0.4=0.25\). Figure~\ref{fig:active_poly_pairing}D shows that this recovery is also achieved with a residual near numerical precision.

\begin{table}[tbp]
\centering
\scriptsize
\input{tables/table_active_polynomial_equations.tex}
\caption{Equation comparison for active polynomial recovery. The table reports the full recovered first-order and second-order polynomial equations. Coefficients removed by thresholding are zero and are omitted.}
\label{tab:active-poly-equations}
\end{table}

%

\subsubsection{Hill growth and transfer}

\begin{figure}[t]
\centering
\includegraphics[width=0.95\linewidth]{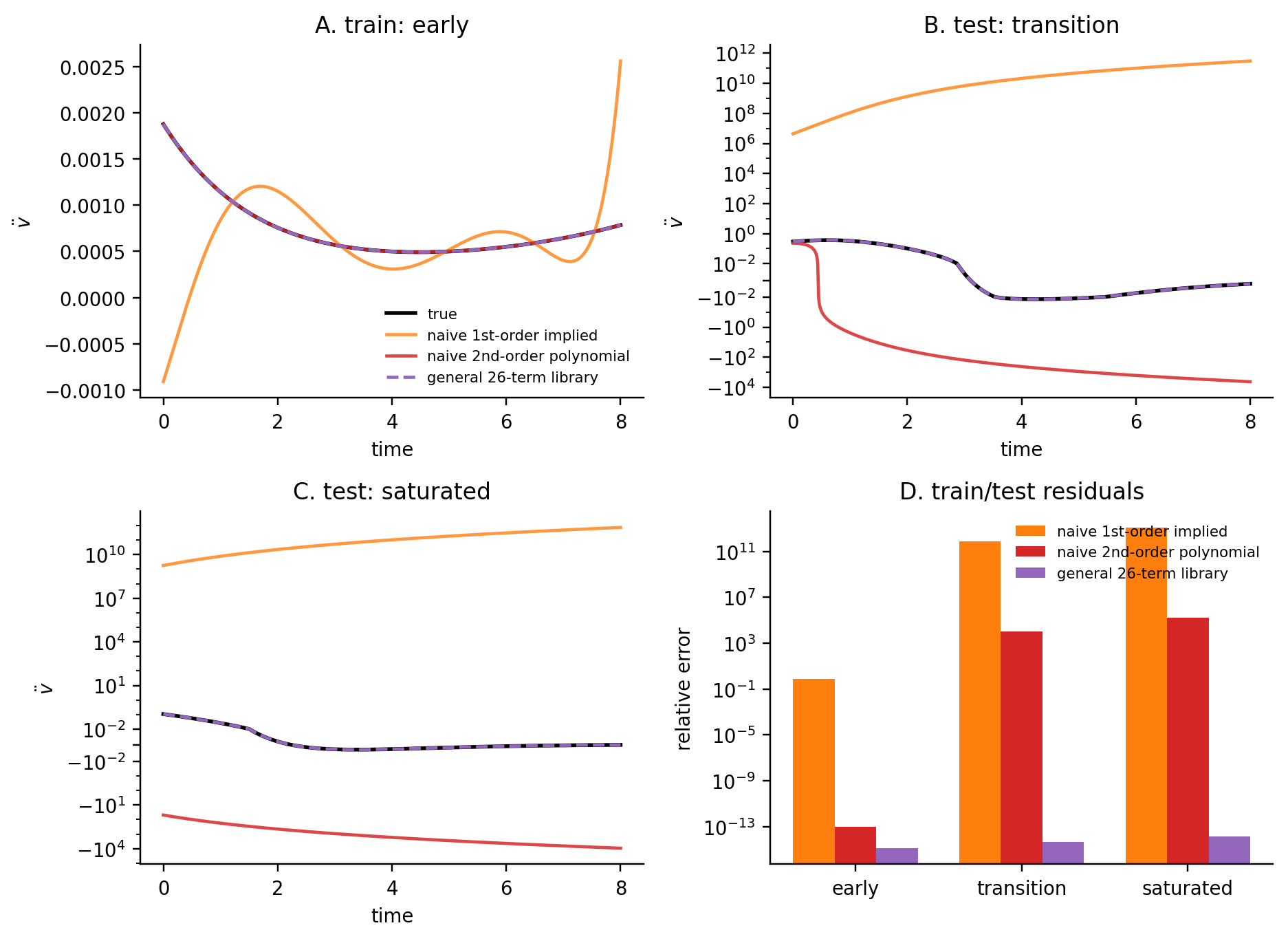}
\caption{Active Hill transfer comparison across early, transition, and saturated regimes. Panels A, B, and C compare recovered accelerations with the truth. Panel D reports relative transfer errors.}
\label{fig:hill_failure}
\end{figure}

The Hill benchmark separates interpolation on one trajectory from recovery of a transferable growth mechanism. The models are trained on the early-growth trajectory, where the Hill law locally resembles a power law. A polynomial surrogate can therefore achieve a small training error without containing the rational saturation structure needed in the transition and saturated regimes. For compact notation, define
\[
G_{n,K}(v)=\frac{v^n}{K^n+v^n}, \qquad
D_{n,K}(v,\dot v)=G'_{n,K}(v)\dot v.
\]
Figure~\ref{fig:hill_failure} should therefore be read as a transfer test. The naive first-order and second-order polynomial models are evaluated using the equations fitted on the early regime. Their errors increase sharply when the trajectory moves through the Hill transition or into saturation. The 26-term library instead selects \(D_{2,1}(v,\dot v)\), \(\dot v\), and \(G_{2,1}(v)\), giving the recovered equation
\[
\ddot v=1.5D_{2,1}(v,\dot v)-0.65\dot v+0.6G_{2,1}(v).
\]
The selected indices identify \(\widehat n=2\) and \(\widehat K=1\). The derivative-coupled coefficient gives \(\widehat a=1.5\), the ratio \(0.6/1.5\) gives \(\widehat\gamma_1=0.4\), and the damping coefficient gives \(\widehat\gamma_2=0.65-0.4=0.25\). Table~\ref{tab:active-hill-equations} compares this equation with the polynomial surrogates. The transfer residual remains near numerical precision in all three regimes, showing that the selected Hill atoms recover the mechanism instead of only an early-regime approximation.

\begin{table}[tbp]
\centering
\scriptsize
\input{tables/table_active_hill_equations.tex}
\caption{Equation comparison for active Hill recovery and transfer.}
\label{tab:active-hill-equations}
\end{table}

\subsection{Aggregate observations}
\subsubsection{Power-law inverse recovery}

For aggregate observations, the measured variable is \(u=v+w\) as opposed to the active state itself. In the power-law case, reconstructing \(v\) from \(\dot u=f(v)\) introduces the fractional powers of \(\dot u\) in \eqref{eq:agg-power}. Polynomial libraries in \(u\) or \((u,\dot u)\) can approximate the observed window, but they do not contain this inverse-dependent structure. The structured aggregate library recovers the fractional-power equation and the associated switching parameters, whereas the polynomial rows in Table~\ref{tab:aggregate_power_equations} should be interpreted as trajectory-level surrogates.

\begin{figure}[tbp]
\centering
\includegraphics[width=0.95\linewidth]{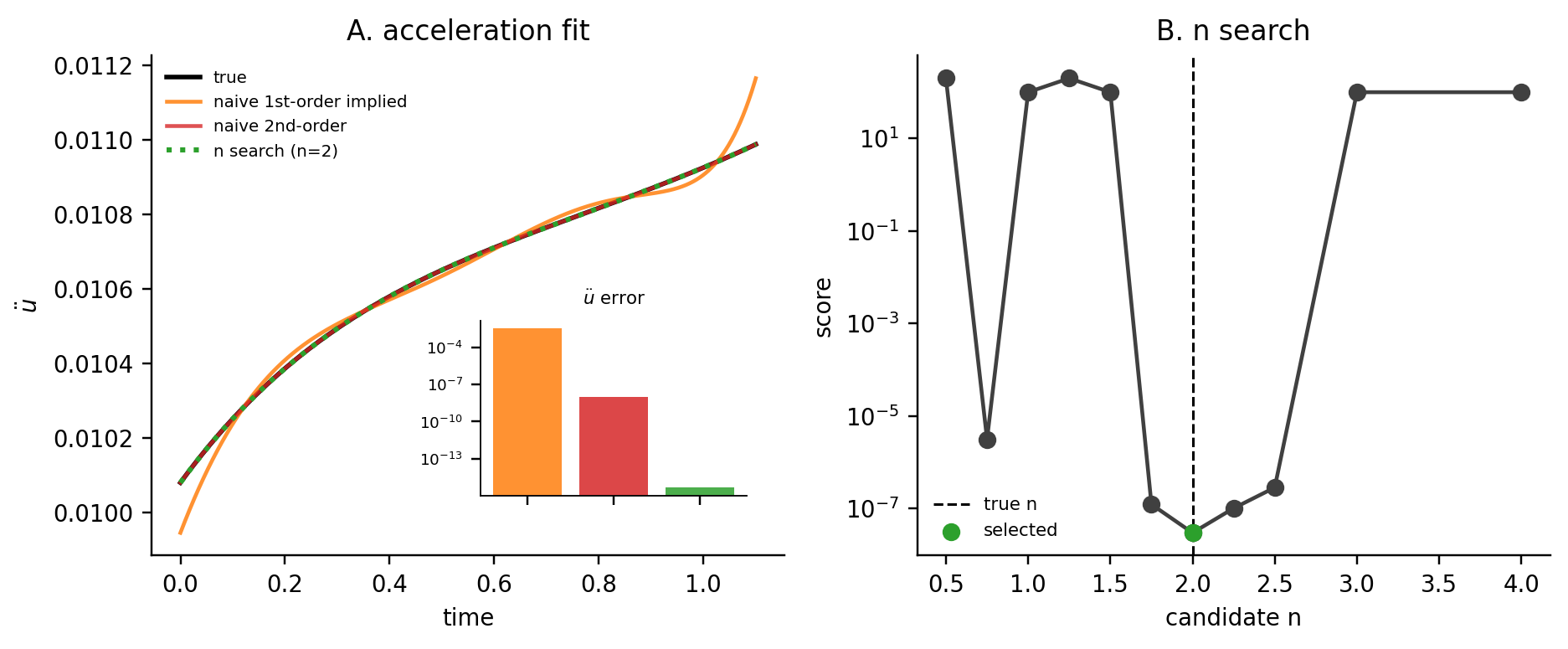}
\caption{Aggregate power-law recovery from \(y=u=v+w\). The correct observable library contains fractional powers of \(\dot u\), whereas blind polynomial libraries in \(u\) or \((u,\dot u)\) are not mechanistically aligned with the aggregate reduction.}
\label{fig:aggregate_power}
\end{figure}

\begin{table}[tbp]
\centering
\scriptsize
\input{tables/table_aggregate_power_equations.tex}
\caption{Equation comparison for aggregate power-law recovery. The rows distinguish the first-order aggregate surrogate, the generic second-order polynomial surrogate, and the structured aggregate \(n\)-search equation.}
\label{tab:aggregate_power_equations}
\end{table}

\subsubsection{Hill nonlinear recovery}

For aggregate Hill growth, the inverse map
\[
    v=K\left(\frac{\dot u}{a-\dot u}\right)^{1/n}
\]
depends on unknown parameters, so there is no fixed linear library to fit first and interpret later. Instead, the problem becomes a nonlinear search with the biological constraint \(0<\dot u<a\). Starting from observations of \(u(t)\), we estimate \(\dot u(t)\) and \(\ddot u(t)\), reconstruct the hidden active state through the inverse Hill map, and then compare the two sides of the reduced observable equation. In other words, the fit is performed by matching the observed left-hand side \(\ddot u(t)\) to the model-predicted right-hand side of the aggregate Hill reduction while searching over \(a\), \(K\), \(n\), \(\gamma_1\), and \(\gamma_2\).

In Appendix~\ref{app:agg-hill} and Figure~\ref{fig:agg_hill}, we show why this is important in practice. The nonlinear search is sensitive to initialization even with synthetic noise-free data. In Table~\ref{tab:aggregate-hill-runs}, one initialization violates the inverse constraint, one gives an admissible but less accurate fit, and the initialization near the true parameters recovers the benchmark values to numerical precision.

Table~\ref{tab:aggregate-hill-runs} separates the residual, parameter error, and admissibility check for each nonlinear aggregate-Hill initialization. This demonstrates when a low residual is accompanied by biologically admissible parameters and when the inverse constraint fails.

\begin{table}[tbp]
\centering
\scriptsize
\resizebox{\textwidth}{!}{\input{tables/table_aggregate_hill_runs.tex}}
\caption{Aggregate-Hill nonlinear search diagnostics for the synthetic benchmark.}
\label{tab:aggregate-hill-runs}
\end{table}

\begin{figure}[h]
\centering
\includegraphics[width=0.90\linewidth]{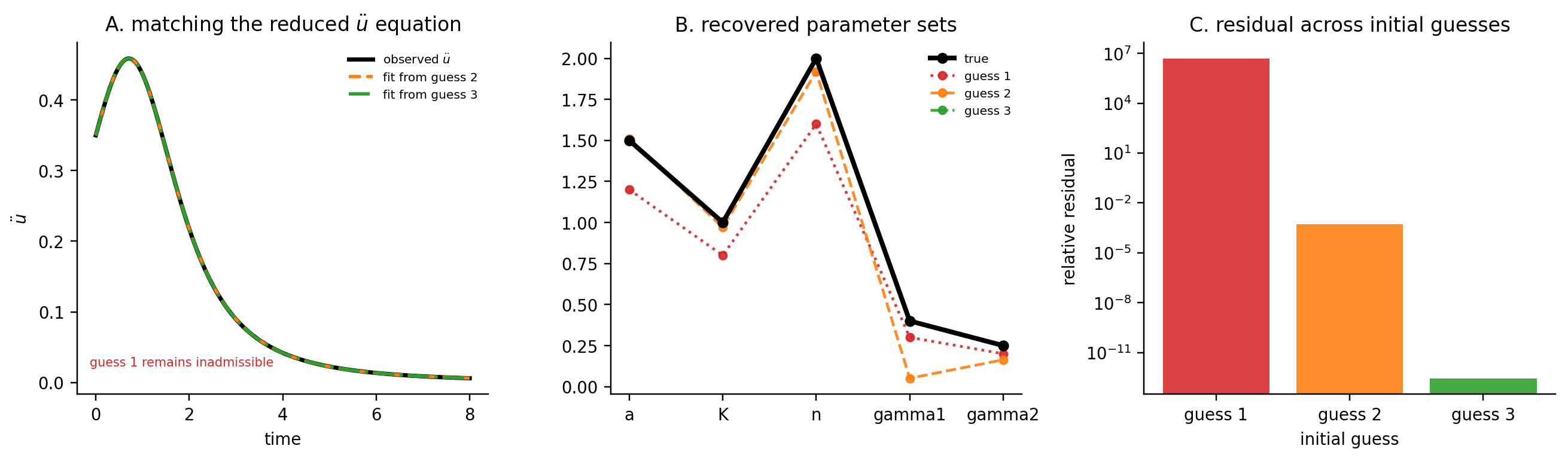}
\caption{Aggregate Hill nonlinear least-squares fits from different initial guesses. A) compares the observed left-hand side \(\ddot u(t)\) with the right-hand side predicted by the recovered aggregate Hill equation. The best fit tracks the observed curve closely, while poorer initial guesses can lead to visibly worse or inadmissible solutions. B) compares the recovered parameter vectors against the true values \(a=1.5\), \(K=1\), \(n=2\), \(\gamma_1=0.4\), and \(\gamma_2=0.25\). C) reports the relative residual for each initial guess, showing the strong sensitivity of the nonlinear search to initialization. The initial guesses were \((1.2, 0.8, 1.6, 0.3, 0.2)\), \((2, 1.5, 2.5, 0.6, 0.5)\), and \((1.55, 1, 2, 0.4, 0.25)\).}
\label{fig:agg_hill}
\end{figure}

\subsubsection{Polynomial branch ambiguity}

The aggregate polynomial case shows a different limitation: even if the model form is correct, the hidden active state may not be uniquely recoverable from the aggregate data. The equation \(f(v)=\dot u\) can have one positive root, several positive roots, or none at all over a given range. When there is more than one admissible positive root, the same measured aggregate growth rate can correspond to different hidden active states.

\begin{figure}[tbp]
\centering
\includegraphics[width=0.95\linewidth]{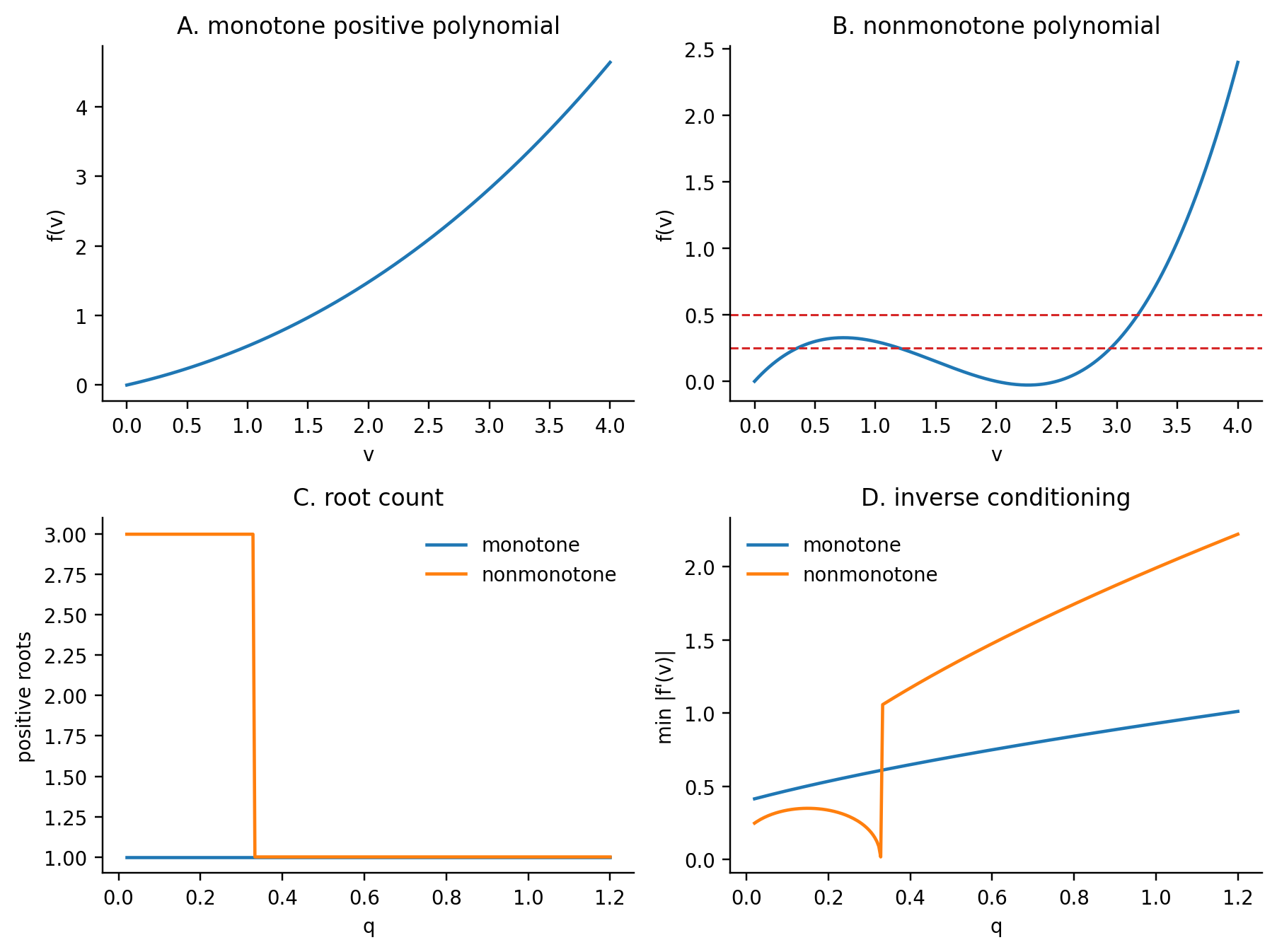}
\caption{Aggregate polynomial branch ambiguity. Panel A shows a monotone polynomial with a unique positive inverse branch. Panel B shows a nonmonotone polynomial for which a single \(q=\dot u\) can correspond to multiple positive active states. Panel C counts positive roots as \(q\) varies. Panel D reports the minimum \(|f'(v)|\) over admissible roots, highlighting regions where inverse reconstruction is ill-conditioned.}
\label{fig:aggregate_poly_branch}
\end{figure}

The aggregate polynomial case gives two pre-regression diagnostics: branch uniqueness and inverse conditioning. In the monotone case shown in Figure~\ref{fig:aggregate_poly_branch}A, each admissible value of \(\dot u\) corresponds to one positive active state, so the inverse branch is unique. In the nonmonotone case shown in Figure~\ref{fig:aggregate_poly_branch}B, the same observed value of \(\dot u\) can correspond to several positive active states. The root count in Figure~\ref{fig:aggregate_poly_branch}C shows this ambiguity directly. A count of one means that the active state can be reconstructed uniquely from aggregate data, while a count greater than one means that aggregate data alone do not determine the hidden active state.

Uniqueness is separate from numerical stability. Even when the positive branch is unique, recovering \(v\) from \(\dot u\) can be unstable if the graph of \(f\) is nearly flat. If a local inverse branch exists, we may write \(v=F(\dot u)\), and the inverse function theorem gives
\[
    \frac{dF}{d\dot u}=\frac{1}{f'(F(\dot u))}.
\]
So, when \(|f'(v)|\) is small, a small perturbation in \(\dot u\) can produce a large change in the reconstructed active state. Figure~\ref{fig:aggregate_poly_branch}D reports the minimum \(|f'(v)|\) along admissible positive roots. Large values indicate stable inversion, while values near zero signal ill conditioning. Aggregate polynomial discovery can therefore fail before sparse regression is attempted, either because the positive inverse branch is nonunique or because the unique branch is too ill-conditioned to use reliably.

\section{Discussion}\label{sec:discussion}

In this work, we show that partial observability changes the equation available for sparse regression before any algorithm is applied. The hidden active-quiescent model is {first-order} in \((v,w)\), but active observation produces a second-order scalar equation and aggregate observation produces an inverse-dependent or branch-dependent scalar equation. This distinction is mathematical, but its practical consequence is biological. The same underlying active-quiescent mechanism can require different data-driven regression libraries depending on whether an experiment measures active markers or aggregate population readouts.

For biological sparse regression, the active-observation results show that an observable-reduction dictionary can cover several growth families without fitting them one at a time. The 26-term library contains polynomial atoms, integer power-law atoms through the polynomial block, and nine Hill atom pairs. It recovers the controlled power-law, polynomial, and Hill observable equations with the expected sparse supports and coefficient maps. The naive polynomial baselines can interpolate restricted trajectories, but they do not provide the same growth-family identification or transfer behavior.

The growth-specific derivations are important because they translate selected terms and coefficients into mechanistic parameters. Sparse regression supplies support selection within the reduction-induced dictionary, while the coefficient maps determine whether the selected support is consistent with active-quiescent switching. Direct parametric fitting remains available when a growth family is known in advance.

Standard SINDy emphasizes sparse selection from a candidate library. Related partial-observation approaches use delay embeddings, learned coordinates, or data assimilation to recover closed-form dynamics when not all state variables are measured \cite{bakarji2023partial,ribera2022model}. The present work is complementary because it derives the scalar observable equation induced by a specified biological observation map before constructing the sparse-regression library. \cite{messenger2021weak,messenger2021weakpde}. SINDy-PI and implicit sparse regression are natural for aggregate Hill and polynomial cases \cite{kaheman2020sindypi,mangan2016biological}, where inverse or implicit relations arise. 

We intentionally use noise-free synthetic trajectories in our present study to separate structural limitations imposed by the biological observation process from statistical limitations caused by finite sampling, measurement noise, and derivative estimation. The resulting analysis therefore represents a best-case biological measurement scenario: failures such as incorrect observable targets, coefficient inconsistency, branch ambiguity, nonidentifiability, and lack of transfer persist even under exact observations and cannot be resolved by adding experimental noise. Conversely, successful recovery in this idealized setting identifies observation protocols and dynamical regimes that are appropriate candidates for subsequent evaluation under realistic noise and sampling conditions. Structural identifiability and observability provide the language \cite{villaverde2019observability,massonis2023distilling} for determining when hidden variables and mechanistic parameters can, in principle, be recovered from output data, whereas weak SINDy and related integral formulations are natural next steps for assessing practical recovery when derivative information must be inferred from noisy or sparsely sampled measurements. Thus, our results suggest that structural observability and branch-admissibility requirements should be established before noise-robust estimation is attempted.

The observable reductions also provide observation-design guidance. Active-state measurements are most favorable when the goal is mechanism recovery because the active observable reduction contains coefficient relations that can be mapped back to switching and growth parameters for selected growth families. Aggregate measurements are common experimentally, but they require reconstructing the hidden active state through \(\dot u=f(v)\). They therefore require branch-admissibility and inverse-conditioning diagnostics in addition to residual fit. Hill-type growth illustrates the role of informative regimes. Early-growth data mainly constrain the lumped ratio \(a/K^n\), while saturated data are weakly informative because \(f(v)\approx a\) and \(f'(v)\approx 0\). Transition-region data near \(v\approx K\) are more informative for separating \(a\), \(K\), and \(n\). Polynomial aggregate reductions show when aggregate data alone are insufficient, because multiple positive roots of \(f(v)=\dot u\) can produce the same observed aggregate growth rate. These examples suggest that sparse active-state markers, multiple initial conditions, replicate trajectories, and perturbations should be viewed as design tools for making observable-reduction-guided sparse regression mechanistically interpretable. This observation-design issue is related to library-conditioning effects in biological equation learning, where restricted trajectory sampling and correlated candidate functions can make sparse recovery unstable even when the correct terms are present in the library \cite{feng2026illconditioning}.

Together, these results motivate a practical workflow for observable-reduction-guided equation learning. One should first specify the observation map and derive the differential equation satisfied by the measured quantity. The candidate library can then be constructed from that reduction, with coefficient relations, branch uniqueness, and inverse conditioning used as structural diagnostics where appropriate. Finally, the recovered model should be evaluated beyond the training regime to distinguish mechanistic recovery from trajectory-specific interpolation.

Our analysis identifies several directions for future research. For active-state observations, the 26-term dictionary relaxes exact growth-family specification, but only over a finite set of candidate functions; broader and adaptively constructed libraries remain to be developed. Extending the same strategy to aggregate observations is more challenging because the corresponding reductions are inverse-dependent or implicit and may require explicit branch selection. Further work should also address PDE models, spatially resolved observations, and stochastic switching to determine how far the dependence of equation learning on the observation map extends beyond the ODE setting. Robustness under realistic sampling and measurement noise should be studied using weak or implicit sparse-regression methods that avoid direct high-order derivative estimation. Finally, resolving branch ambiguity in aggregate polynomial models will require either additional measurements, such as sparse active-state markers, or experimentally justified constraints on the admissible growth law.

\section{Conclusion}\label{sec:conclusion}

For partially observed active-quiescent systems, the hidden two-state dynamics and the differential equation satisfied by the measured signal are generally not the same. Eliminating the unobserved compartment changes both the appropriate regression target and the structure of the candidate library. Under active-state observation, the resulting second-order scalar equation contains coefficient relations that can be used to identify growth-law structure, recover switching parameters, and test mechanistic consistency. Under aggregate observation, the reduction is more restrictive because the hidden active state must be reconstructed through the growth law, introducing inverse dependence, branch ambiguity, and conditioning constraints. These results show that a sparse-regression model should not be interpreted from its residual or selected terms. Mechanistic interpretation requires first incorporating the biological observation map into the regression target, library construction, and validation diagnostics.
\appendix
\section{Derivations of the Observable Reductions}\label{app:observable-derivations}

This appendix gives the algebraic steps leading to the active and aggregate observable reductions used in Section~\ref{sec:reductions}.

\subsection{Active-state reduction}

Assume \(\gamone>0\) and \(y=v\). Solving \eqref{eq:model-v} for the hidden quiescent state gives
\begin{equation}\label{eq:app-w-from-v}
    w=\frac{\vdot-f(v)+\gamtwo v}{\gamone}.
\end{equation}
Differentiating \eqref{eq:model-v} gives
\begin{equation}\label{eq:app-vddot-step1}
    \vddot=f'(v)\vdot-\gamtwo\vdot+\gamone\dot w.
\end{equation}
Using \eqref{eq:model-w}, we have
\begin{equation}\label{eq:app-vddot-step2}
    \gamone\dot w=-\gamone^2w+\gamone\gamtwo v.
\end{equation}
Substituting \eqref{eq:app-w-from-v} into \eqref{eq:app-vddot-step2} gives
\begin{equation}\label{eq:app-vddot-step3}
    \gamone\dot w=-\gamone\vdot+\gamone f(v)-\gamone\gamtwo v+\gamone\gamtwo v.
\end{equation}
Therefore
\begin{equation}\label{eq:app-active-general}
    \boxed{
    \vddot=\bigl(f'(v)-\gamone-\gamtwo\bigr)\vdot+\gamone f(v).
    }
\end{equation}
which is \eqref{eq:active-general}.

\subsection{Aggregate-state reduction}

Let \(u=v+w\). Adding \eqref{eq:model-v} and \eqref{eq:model-w} gives
\begin{equation}\label{eq:app-u-dot}
    \udot=f(v).
\end{equation}
Since \(w=u-v\), the active equation becomes
\begin{equation}\label{eq:app-vdot-aggregate}
    \vdot=f(v)-\gamtwo v+\gamone(u-v).
\end{equation}
Using \(\udot=f(v)\), this can be written as
\begin{equation}\label{eq:app-vdot-aggregate2}
    \vdot=\udot+\gamone u-(\gamone+\gamtwo)v.
\end{equation}
Differentiating \eqref{eq:app-u-dot} gives
\begin{equation}\label{eq:app-uddot-step}
    \uddot=f'(v)\vdot.
\end{equation}
Substituting \eqref{eq:app-vdot-aggregate2} into \eqref{eq:app-uddot-step} gives the differential algebraic aggregate reduction
\begin{subequations}\label{eq:app-agg-dae}
\begin{align}
    \udot &= f(v),\\
    \uddot &= f'(v)\left[\udot+\gamone u-(\gamone+\gamtwo)v\right].
\end{align}
\end{subequations}
If \(f\) has a locally admissible inverse branch \(F\), then \(v=F(\udot)\). Substitution into \eqref{eq:app-agg-dae} gives
\begin{equation}\label{eq:app-agg-general}
    \uddot=f'\!\left(F(\udot)\right)
    \left[\udot+\gamone u-(\gamone+\gamtwo)F(\udot)\right],
\end{equation}
which is \eqref{eq:agg-general}.

\section{Derivations for Observation-Specific Libraries}\label{app:growth-library-derivations}

This appendix records the substitutions used to obtain the growth law specific libraries in Section~\ref{sec:libraries}. The main text keeps the final targets and coefficient maps because those are the objects used for sparse regression.

\subsection{Power-law growth}

For \(f(v)=av^n\), we have \(f'(v)=anv^{n-1}\). Substitution into the active reduction \eqref{eq:active-general} gives
\begin{equation}\label{eq:app-active-power}
    \boxed{
    \vddot=anv^{n-1}\vdot-(\gamone+\gamtwo)\vdot+\gamone av^n.
    }
\end{equation}
This is the target \eqref{eq:active-power}. Matching coefficients in
\[
    \vddot=c_1v^{n-1}\vdot+c_2\vdot+c_3v^n
\]
gives \(c_1=an\), \(c_2=-(\gamone+\gamtwo)\), and \(c_3=\gamone a\), which yields \eqref{eq:active-power-map}.

For aggregate observation, \(\udot=av^n\) gives
\begin{equation}\label{eq:app-power-inverse}
    v=\left(\frac{\udot}{a}\right)^{1/n}.
\end{equation}
Also,
\begin{equation}\label{eq:app-power-fprime}
    f'(v)=anv^{n-1}=n a^{1/n} \udot^{(n-1)/n}.
\end{equation}
Substituting \eqref{eq:app-power-inverse} and \eqref{eq:app-power-fprime} into \eqref{eq:agg-general} gives
\begin{align}\label{eq:app-agg-power-expanded}
    \uddot
    &=n a^{1/n}\udot^{(n-1)/n}\left[\udot+\gamone u-(\gamone+\gamtwo)\left(\frac{\udot}{a}\right)^{1/n}\right]\\
    &=n a^{1/n}\udot^{(2n-1)/n}+n\gamone a^{1/n}u \udot^{(n-1)/n}-n(\gamone+\gamtwo)\udot.
\end{align}
This is \eqref{eq:agg-power}. Matching coefficients gives \(c_1=na^{1/n}\), \(c_2=n\gamone a^{1/n}\), and \(c_3=-n(\gamone+\gamtwo)\), which yields \eqref{eq:agg-power-map}.

\subsection{Hill-type growth}

For
\[
    f(v)=\frac{av^n}{K^n+v^n},
\]
the derivative is
\begin{equation}\label{eq:app-hill-derivative}
    f'(v)=\frac{anK^n v^{n-1}}{(K^n+v^n)^2}.
\end{equation}
Substituting \eqref{eq:app-hill-derivative} into \eqref{eq:active-general} gives the active Hill target \eqref{eq:active-hill}. For fixed \((n,K)\), the terms in \eqref{eq:active-hill} split into the atoms \(D_{n,K}(v,\dot v)\), \(\dot v\), and \(G_{n,K}(v)\) defined in \eqref{eq:active-hill-atoms}. Coefficient matching gives \(c_D=an\), \(c_G=\gamone a\), and \(c_{\dot v}=-(\gamone+\gamtwo)\).

For aggregate observation, solving
\[
    \udot=\frac{av^n}{K^n+v^n}
\]
for \(v>0\) gives
\begin{equation}\label{eq:app-hill-inverse}
    v=K\left(\frac{\udot}{a-\udot}\right)^{1/n},
    \qquad 0<\udot<a.
\end{equation}
Using \(\udot=f(v)\), the derivative can also be written as
\begin{equation}\label{eq:app-hill-fprime-q}
    f'(v)=\frac{n\udot(a-\udot)}{av}.
\end{equation}
Substitution of \eqref{eq:app-hill-inverse} and \eqref{eq:app-hill-fprime-q} into \eqref{eq:agg-general} gives the aggregate Hill target \eqref{eq:agg-hill}. Because \(a\), \(K\), and \(n\) enter the inverse reconstruction, this target is nonlinear in the unknown parameters.

\subsection{Polynomial growth}

For
\[
    f(v)=\sum_{j=m}^{p}a_jv^j,
\]
we have
\begin{equation}\label{eq:app-poly-derivative}
    f'(v)=\sum_{j=m}^{p}j a_jv^{j-1}.
\end{equation}
Substitution into \eqref{eq:active-general} gives
\begin{equation}\label{eq:app-active-poly}
    \vddot=
    \sum_{j=m}^{p}j a_jv^{j-1}\vdot
    -(\gamone+\gamtwo)\vdot
    +\gamone\sum_{j=m}^{p}a_jv^j.
\end{equation}
which is \eqref{eq:active-poly}. The coefficient of \(v^{j-1}\dot v\) is \(j a_j\), and the coefficient of \(v^j\) is \(\gamone a_j\). Eliminating \(a_j\) gives the {degree-wise} estimate \eqref{eq:gamma1-pair}.

For aggregate observation, substituting the polynomial law into \eqref{eq:agg-dae} gives
\begin{subequations}\label{eq:app-agg-poly-dae}
\begin{align}
    \udot &=\sum_{j=m}^{p}a_jv^j,\\
    \uddot
    &=
    \left(\sum_{j=m}^{p}j a_jv^{j-1}\right)
    \left[\udot+\gamone u-(\gamone+\gamtwo)v\right].
\end{align}
\end{subequations}
Recovering the hidden active state requires solving
\begin{equation}\label{eq:app-poly-root}
    \sum_{j=m}^{p}a_jv^j-\udot=0.
\end{equation}
The positive branch conditions and conditioning of this reconstruction are discussed in Appendix~\ref{app:poly-branches}.

\section{Aggregate Hill Nonlinear Search}\label{app:agg-hill}

The aggregate Hill case is not an ordinary linear SINDy problem because the inverse
\[
    v=K\left(\frac{\dot u}{a-\dot u}\right)^{1/n}
\]
depends on the unknown parameter \(a\). For a candidate parameter vector
\[
    \theta=(a,K,n,\gamma_1,\gamma_2),
\]
the strong-form residual is
\begin{equation}\label{eq:agg-hill-resid}
    R(t;\theta)
    =
    \ddot u(t)
    -
    \frac{n\dot u(t)\bigl(a-\dot u(t)\bigr)}{a\,v(t;\theta)}
    \left[
    \dot u(t)+\gamma_1u(t)-(\gamma_1+\gamma_2)v(t;\theta)
    \right],
\end{equation}
with \(v(t;\theta)=K(\dot u(t)/(a-\dot u(t)))^{1/n}\). Invalid parameter values are penalized when \(a\leq 0\), \(K\leq0\), \(n\leq0\), \(\gamma_i<0\), or when \(0<\dot u(t)<a\) fails. The nonlinear least-squares problem can be sensitive to initial guesses because \(a\) controls both the inverse and the prefactor in \eqref{eq:agg-hill-resid}.

\section{Positive Branches for Aggregate Polynomial Reductions}\label{app:poly-branches}

For aggregate observations, the first reduced equation is
\[
    \dot u=f(v).
\]
When \(f\) is polynomial,
\[
    f(v)=\sum_{j=m}^{p}a_jv^j,
\]
recovering the hidden active state requires solving
\[
    \sum_{j=m}^{p}a_jv^j-\dot u=0.
\]
Although this degree \(p\) equation may have complex or negative roots, biological admissibility requires \(v>0\) for positive population densities. Therefore the relevant root set is
\[
    \mathcal A(\dot u)=\{v>0:f(v)=\dot u\}.
\]

If \(|\mathcal A(\dot u)|=1\) over the observed range, then the aggregate data determine a unique positive branch \(v=\phi(\dot u)\). If \(|\mathcal A(\dot u)|\geq2\), there is biological branch ambiguity. At least two positive active states generate the same aggregate growth rate. If \(|\mathcal A(\dot u)|=0\), the observation is inconsistent with the assumed polynomial law on the positive domain, possibly because of model mismatch or derivative-estimation error.

Local regularity is separate from uniqueness. If \(f'(\phi(\dot u))\neq0\), then the inverse function theorem gives a smooth local branch and
\[
    \frac{d\phi}{d\dot u}=\frac{1}{f'(\phi(\dot u))}.
\]
Thus reconstruction is {ill-conditioned} when \(f'\) is close to zero. For aggregate polynomial sparse regression, {branch reconstruction} requires a unique positive branch, while stable reconstruction requires that branch to remain away from critical points. Thus positive-branch uniqueness is a reconstruction condition, while stable use of that branch also requires the inverse to be well conditioned over the observed range.

\section{Mathematical Background: Hidden States, Observable Reduction, and Sparse Regression}
\label{app:mathematical_background}

\subsection{Hidden states and observable reduction}

Many biological systems contain distinct subpopulations whose internal states cannot all be measured experimentally. In the active-quiescent model considered here, the full state is
\[
x(t)=
\begin{pmatrix}
v(t)\\
w(t)
\end{pmatrix},
\]
where \(v(t)\) denotes the actively proliferating population and \(w(t)\) denotes the quiescent population. Although cells may switch between these states, growth occurs only through the active compartment. An experiment generally records an observable
\[
y(t)=P(x(t)),
\]
instead of the complete state \(x(t)\). For example, a proliferation marker may approximately measure only the active population,
\[
y(t)=v(t),
\]
whereas total cell number, optical density, tumor burden, or aggregate fluorescence may measure the combined population
\[
y(t)=u(t)=v(t)+w(t).
\]
Thus, partial observation measures one component of the system, while aggregate observation combines multiple biological states into a single measured quantity.

Differential equations provide a mechanistic description of how a biological system changes over time by relating its current state to its rate of change. A differential equation relates the time derivative of the observable, such as \(\dot{y}(t)\), to terms on the right hand side that represent the biological mechanisms influencing how the system's variables change over time. The unobserved state can still influence the evolution of the observable, so the measured quantity generally does not satisfy the same differential equation as the full hidden system. \emph{Observable reduction} is the process of eliminating the hidden variables to derive an equation involving only the measured quantity and its time derivatives. For active-state observation, eliminating \(w(t)\) typically produces a differential equation with second order time derivatives of the form
\[
\ddot{v}=g_{\mathrm{act}}(v,\dot{v}).
\]
For aggregate observation, the identity
\[
\dot{u}=f(v)
\]
shows that the hidden active population must instead be reconstructed from the observed aggregate growth rate. When an admissible inverse branch \(F\) exists, this reconstruction takes the form
\[
v=F(\dot{u}),
\]
leading to an aggregate observable equation
\[
\ddot{u}=g_{\mathrm{agg}}(u,\dot{u}).
\]
The form of the reduced equation therefore depends on what the experiment measures, even when the underlying biological mechanism is unchanged.

\noindent\emph{For further reading on hidden-state models, partial observations, observable reduction, and identifiability, see
\cite{somacal2022hidden,lu2022partial,bakarji2023partial,ribera2022model,villaverde2019observability,massonis2023distilling}.}

\subsection{Equation learning and sparse regression}

The goal of equation learning is to infer a mathematical model directly from time-series data. Instead of assuming a single equation in advance, the method begins with a collection of plausible mathematical terms, called a \emph{library}, and determines a parsimonious subset that best explains the observed dynamics. For a measured variable \(y(t)\), a candidate library might be written as
\[
\Theta(y)=
\begin{bmatrix}
1 & y & y^2 & y^3 & \cdots
\end{bmatrix}.
\]
More specialized libraries may include interaction terms, derivatives, fractional powers, or saturating functions suggested by the underlying biology. The goal of equation learning is automated model selection, i.e., to identify a model that uses as few terms as possible among the library to reproduce the measured dynamics.

Sparse Identification of Nonlinear Dynamics, or SINDy, expresses this idea as a regression problem. Given measurements \(y(t_i)\) and an estimated regression target \(z(t_i)\), the candidate functions are evaluated at the sampled data to form a matrix \(\Theta\), and SINDy seeks a sparse coefficient vector \(\xi\) satisfying
\[
z\approx\Theta\xi.
\]
In standard first-order SINDy, \(z=\dot{y}\). For a partially observed system, however, observable reduction may show that the appropriate target is instead \(z=\ddot{y}\), with a library constructed from \(y\) and \(\dot{y}\):
\[
\ddot{y}\approx\Theta(y,\dot{y})\xi.
\]
Sparsity means that most entries of \(\xi\) are zero, leaving only a small number of selected terms. A low regression error alone does not guarantee that the selected equation represents the biological mechanism because an overly flexible or incorrectly constructed library may reproduce one trajectory while failing to recover meaningful parameters or predict a new experimental regime. The regression target and candidate library should therefore be derived from the observable available in the experiment.

\noindent\emph{For further reading on SINDy, sparse equation learning, implicit formulations, and biological applications, see
\cite{brunton2016sindy,mangan2016biological,kaheman2020sindypi,messenger2021weak,kaptanoglu2022pysindy}.}

\section*{Declarations}

\paragraph{Authorship and contributorship}
KBF conceived the study. KCN and KBF designed the research. KCN performed the analysis and generated the results. KCN and KBF interpreted the results. KBF procured
funding for the research. KCN drafted the initial manuscript, and all authors reviewed, edited, and approved
the final manuscript.

\paragraph{Conflicts of interest}
The authors declare that they have no conflict of interest.

\paragraph{Data and code availability}
The data and code used to support the findings of this study is available in the project GitHub repository: \\
https://github.com/kcnguyen3191/quiescent-eql.git

\paragraph{Funding}
KBF was funded in part by the National Science Foundation under grant numbers DMS-2327836, DMS-2342344, and DMS-2424748, and by the National Institute of Allergy and Infectious Diseases of the National Institutes of Health under grant number 1U54AI191253-01.

\bibliographystyle{siamplain}
\bibliography{references}

\end{document}

%% file: tables/table_active_power_equations.tex
\begin{tabular}{p{0.22\linewidth}p{0.70\linewidth}}
\toprule
Model & Equation \\
\midrule
True observable & $\begin{aligned}\ddot v &= 1.6\,v\dot v - 0.65\,\dot v + 0.32\,v^2\end{aligned}$ \\
\midrule
Naive first order & $\begin{aligned}\dot v &= 0.122 - 0.799\,v + 3.43\,v^2 - 22.4\,v^3 + 106\,v^4 - 160\,v^5\end{aligned}$ \\
\midrule
Naive second order & $\begin{aligned}\ddot v &= 1.6\,v\dot v - 0.65\,\dot v + 0.32\,v^2\end{aligned}$ \\
\midrule
General 26-term library & $\begin{aligned}\ddot v &= 1.6\,v\dot v - 0.65\,\dot v + 0.32\,v^2\end{aligned}$ \\
\bottomrule
\end{tabular}

%% file: tables/table_active_polynomial_equations.tex
\begin{tabular}{p{0.28\linewidth}p{0.64\linewidth}}
\toprule
Model & Equation \\
\midrule
True observable & $\begin{aligned}\ddot v &= 0.16\,v\dot v - 0.03\,v^2\dot v - 0.65\,\dot v + 0.032\,v^2 - 0.004\,v^3\end{aligned}$ \\
\midrule
Naive first-order polynomial library & $\begin{aligned}\dot v &= 6.39 - 29.8\,v + 51.5\,v^2 - 39.8\,v^3 + 12.4\,v^4 - 0.708\,v^5\end{aligned}$ \\
\midrule
Naive second-order polynomial library & $\begin{aligned}\ddot v &= 0.16\,v\dot v - 0.03\,v^2\dot v - 0.65\,\dot v + 0.032\,v^2 - 0.004\,v^3\end{aligned}$ \\
\midrule
General 26-term library & $\begin{aligned}\ddot v &= 0.16\,v\dot v - 0.03\,v^2\dot v - 0.65\,\dot v + 0.032\,v^2 - 0.004\,v^3\end{aligned}$ \\
\bottomrule
\end{tabular}

%% file: tables/table_active_hill_equations.tex
\begin{tabular}{p{0.17\linewidth}p{0.57\linewidth}p{0.16\linewidth}}
\toprule
Model & Equation & Diagnostic \\
\midrule
True observable & $\begin{aligned}\ddot v={}&1.5D_{2,1}(v,\dot v)-0.65\dot v+0.6G_{2,1}(v)\end{aligned}$ & \(a=1.5,\ K=1,\ n=2\) \\
\midrule
Naive first order & $\begin{aligned}\dot v &= -17.4 + 1.429\times 10^{3}\,v - 4.680\times 10^{4}\,v^2\\ &\quad + 7.643\times 10^{5}\,v^3 - 6.224\times 10^{6}\,v^4 + 2.022\times 10^{7}\,v^5\end{aligned}$ & Degree 5 polynomial surrogate \\
\midrule
Naive second order & $\begin{aligned}\ddot v &= -0.65\,\dot v + 0.6\,v^2 + 3\,v\dot v\\ &\quad + 0.00185\,\dot v^2 - 4.352\times 10^{-4}\,v^3 - 2.865\times 10^{-4}\,v^2\dot v\\ &\quad - 0.127\,v\dot v^2 - 0.102\,\dot v^3 - 0.592\,v^4\\ &\quad - 6\,v^3\dot v + 1.77\,v^2\dot v^2 + 3.67\,v\dot v^3\\ &\quad + 0.645\,\dot v^4\end{aligned}$ & 13-term polynomial surrogate \\
\midrule
General 26-term library & $\begin{aligned}\ddot v={}&1.5D_{2,1}(v,\dot v)-0.65\dot v+0.6G_{2,1}(v)\end{aligned}$ & Selected \(G_{2,1},D_{2,1},\dot v\) \\
\bottomrule
\end{tabular}

%% file: tables/table_aggregate_power_equations.tex
\begin{tabular}{p{0.22\linewidth}p{0.70\linewidth}}
\toprule
Model & Equation \\
\midrule
True observable & $\begin{aligned}\ddot u &= 1.789\,(\dot u)^{3/2} + 0.716\,u(\dot u)^{1/2} - 1.3\,\dot u\end{aligned}$ \\
\midrule
Naive first order & $\begin{aligned}\dot u &= - 1.22\times 10^{4} + 1.96\times 10^{5}\,u - 1.26\times 10^{6}\,u^{2}\\&\quad + 4.07\times 10^{6}\,u^{3} - 6.56\times 10^{6}\,u^{4} + 4.23\times 10^{6}\,u^{5}\end{aligned}$ \\
\midrule
Naive second order & $\begin{aligned}\ddot u &= - 2.997 + 123.847\,u - 650.972\,\dot u\\&\quad - 645.986\,u^{2} + 4.87\times 10^{3}\,u \dot u - 1.54\times 10^{3}\,(\dot u)^{2}\\&\quad + 884.961\,u^{3} - 8.96\times 10^{3}\,u^{2} \dot u + 8.73\times 10^{3}\,u (\dot u)^{2}\\&\quad + 1.37\times 10^{4}\,(\dot u)^{3}\end{aligned}$ \\
\midrule
$n$-search SINDy ($n=2$) & $\begin{aligned}\ddot u &= 1.789\,(\dot u)^{3/2} + 0.716\,u(\dot u)^{1/2} - 1.3\,\dot u\end{aligned}$ \\
\bottomrule
\end{tabular}

%% file: tables/table_aggregate_hill_runs.tex
\begin{tabular}{lllll}
\toprule
Initial guess & Residual \(R\) & Parameter error \(E_{\mathrm{param}}\) & Admissible \(0<\dot u<a\) & Comment \\
\midrule
\((1.2,0.8,1.6,0.3,0.2)\) & \(4.963\times 10^{6}\) & \(2.012\times 10^{-1}\) & no & inverse constraint fails \\
\midrule
\((2,1.5,2.5,0.6,0.5)\) & \(5.181\times 10^{-4}\) & \(1.354\times 10^{-1}\) & yes & admissible fit \\
\midrule
\((1.55,1,2,0.4,0.25)\) & \(2.810\times 10^{-13}\) & \(1.402\times 10^{-12}\) & yes & admissible fit \\
\bottomrule
\end{tabular}